\documentclass[sn-mathphys-num]{sn-jnl}

\usepackage[T1]{fontenc}
\usepackage[utf8]{inputenc}
\usepackage[english]{babel}
\usepackage{comment}
\usepackage{tabularx}
\usepackage{amsmath,amssymb,amsfonts}%
\usepackage{amsthm}%
\usepackage{mathrsfs}%
\usepackage[title]{appendix}%
\usepackage{xcolor}%
\usepackage{textcomp}%
\usepackage{manyfoot}%
\usepackage{hyperref,float}
\usepackage{bbold} 
\usepackage{hyperref,float,multirow}
\usepackage{colortbl}
\usepackage{array}
\usepackage{multirow}
\usepackage{booktabs}

\usepackage{makeidx,geometry}
\usepackage{listings}%
\usepackage{algorithm,algorithmic}

\newcommand{\BEGIN}{ \STATE \fbox{\textbf{procedure}}}

\newcommand{\Input}{\STATE \textbf{Input}:}
\newcommand{\Output}{\STATE \textbf{Output: }}

\theoremstyle{thmstyleone}%

\theoremstyle{thmstyletwo}%

\theoremstyle{thmstylethree}%

\begin{document}
	\title[Article Title]{Penalized estimation of GEV parameters for extreme quantile regression}

	\author*[1]{\fnm{Lucien M. } \sur{Vidagbandji}}\email{mahutin-lucien.vidagbandji@univ-lehavre.fr}
	\author[1]{\fnm{Alexandre} \sur{Berred}}\email{alexandre.berred@univ-lehavre.fr}
	
	\author[2]{\fnm{Cyrille} \sur{Bertelle}}\email{cyrille.bertelle@univ-lehavre.fr}
	\author[2]{\fnm{Laurent} \sur{Amanton}}\email{laurent.amanton@univ-lehavre.fr}
%
	\affil[1]{\orgdiv{Université Le Havre Normandie, Normandie Univ.}, \orgname{LMAH UR 3821}, \orgaddress{ \city{76600},  \state{Le Havre}, \country{France}}}
	
	\affil[2]{\orgdiv{Université Le Havre Normandie, Normandie Univ.}, \orgname{LITIS},\orgaddress{ \city{76600},  \state{Le Havre}, \country{France}}}
	\abstract{ 
Quantile regression (QR) relies on the estimation of conditional quantiles and explores the relationships between independent and dependent variables. At high probability levels, classical QR methods face extrapolation difficulties due to the scarcity of data in the tail of the distribution. Another challenge arises when the number of predictors is large and the quantile function exhibits a complex structure. In this work, we propose an estimation method designed to overcome these challenges. To enhance extrapolation in the tail of the conditional response distribution, we model block maxima using the generalized extreme value (GEV) distribution, where the parameters depend on covariates. To address the second challenge, we adopt an approach based on generalized random forests (grf) to estimate these parameters. Specifically, we maximize a penalized likelihood, weighted by the weights obtained through the grf method. This penalization helps overcome the limitations of the maximum likelihood estimator (MLE) in small samples, while preserving its optimality in large samples. The effectiveness of our method is validated through comparisons with other approaches in simulation studies and an application to U.S. wage data.
}


	\maketitle
	\section{Introduction}\label{sec1}
	
The study of extreme phenomena has become an unavoidable necessity in the current context, where their impacts are increasingly concerning. Whether related to climate change \cite{arnell_unbiased_1988}, financial crises, or technological  failures, extreme events—though rare by definition—often lead to severe and disproportionate consequences. Positioned in the extreme tails of statistical distributions, these events exhibit unique characteristics that cannot be captured by average or typical behaviors. Their modeling and understanding, therefore, require specialized tools capable of grasping the complexity and dynamics inherent to these phenomena.  Classical statistical methods, which primarily focus on analyzing the central values of distributions, provide powerful frameworks for examining average trends and central relationships between variables. However, they prove inadequate when it comes to exploring or predicting behaviors in the extreme regions of distributions. The tails of distributions, where extreme values reside, are often underrepresented in standard models due to the low data density and the difficulty of capturing the complex relationships that prevail in these regions.   In response to these limitations, new methodologies have emerged to address the specific needs of extreme event analysis. Among them, extreme quantile regression stands out for its relevance and constitutes the main focus of this work.

 The quantile regression model, as introduced by \cite{koenker_regression_1978}, aims to estimate a conditional quantile from a sample \((X_1, Y_1), \ldots, (X_n, Y_n)\), representing independent copies of the random vector \((X, Y)\), where \(X \in \mathcal{X} \subset \mathbb{R}^p\) and \(Y \in \mathcal{Y} \subset \mathbb{R}\). This model allows for the estimation of the conditional quantile \(\mathcal{Q}_x(\tau)\), defined as  
\begin{equation}\label{quant_inverseF}
		\mathcal{Q}_{x}(\tau)= F_{Y|X=x}^{-1}(\tau), 
	\end{equation}
where \(x \in \mathcal{X}\), \(\tau \in (0,1)\), and \(F_{Y|X=x}^{-1}\) represents the generalized inverse of the conditional cumulative distribution function of \(Y\), given \(X = x\). This model has revolutionized the statistical approach by providing a powerful alternative to classical regression \cite{beyerlein_quantile_2014}.  Unlike the latter, which is limited to estimating the conditional mean \(E(Y|X=x)\), quantile regression allows for the estimation of any conditional quantile. This flexibility enables a deeper understanding of heterogeneous relationships between variables, highlighting dynamics that traditional methods fail to capture. Numerous quantile regression models have been proposed in the literature, including those by  \cite{angrist_quantile_2006}, \cite{koenker_quantile_2017}, \cite{benziadi_recursive_2016}, \cite{cade_gentle_2003}, and \cite{wang_tail_2009}. These approaches are effective for moderate quantiles, where sufficient observations are available in the sample. However, when focusing on extreme quantiles, located in the tails of the distribution, new challenges arise.

Classical quantile regression methods perform well when $\tau_n \to 1$ and $n(1-\tau_n)\to \infty$ as $n \to \infty$. In contrast, when $\tau_n \to 1$ and $n(1-\tau_n) \to d \in (0,+\infty)$, estimation requires extrapolation into the tails. This is problematic because observations are scarce in these regions and the statistical properties of extremes differ fundamentally from those in the central part of the distribution. Extreme quantile regression, leveraging advancements from extreme value theory (EVT) and incorporating robust estimation techniques, emerges as an appropriate response to these challenges. Several studies have explored the integration of extreme value theory in this framework, including those by \cite{chernozhukov_extremal_2005}, \cite{zheng_composite_2017}, \cite{zhu_extreme_2022}, \cite{schaumburg_predicting_2012}, \cite{saulo_new_2022}, \cite{chernozhukov_fast_2020}, and \cite{kithinji_adjusted_2021}. A more detailed presentation of extreme quantile regression approaches can be found in the work of \cite{chernozhukov_extremal_2017}.

Even when sufficient data is available, classical regression methods can encounter difficulties if the conditional quantile function is nonlinear or heterogeneous. In such cases, overly simplistic models risk introducing bias \cite{gnecco_extremal_2024}. To address this issue,
several models combining statistical learning methods have been proposed.
Among the notable approaches are those based on neural networks, such as \cite{cannon_non-crossing_2018}, as well as methods relying on decision trees and random forests, including those proposed by \cite{chaudhuri_nonparametric_2002}, \cite{meinshausen_quantile_2006}, \cite{buchinsky_recent_1998}, and \cite{athey_generalized_2019}. Additionally, studies utilizing other machine learning techniques, such as those by \cite{ye_non-parametric_nodate}, \cite{yao_spatiotemporal_2022}, and \cite{tyralis_hydrological_2019}, also provide significant contributions.

In high-dimensional settings, identifying neighbors close to $x$ becomes problematic, and kernel or nearest-neighbor methods suffer from the curse of dimensionality. Recently, several approaches have combined EVT with statistical learning to simultaneously address extrapolation in the tails, complexity of quantile fonction, and high dimensionality of covariate space. Examples include generalized additive models proposed by \cite{youngman_generalized_2019}, gradient boosting introduce by \cite{velthoen_gradient_2023}, neural networks  by \cite{pasche_neural_2023}, and generalized random forests proposed by \cite{gnecco_extremal_2024} and \cite{Vidagbandji2025}.  To enable extrapolation in the distribution tails, these authors use the generalized Pareto distribution (GPD) within the Peak-over-Threshold (POT) approach of extreme value theory. Specifically, these models approximate the conditional distribution \(Y|X=x\) in expression (\ref{quant_inverseF}) with the GPD, where the parameters depend on the covariate \(x\), and thus obtain the extreme conditional quantile using various techniques for estimating the GPD parameters.

In this paper, we adopt the block maxima approach from EVT (see \cite{de_haan_extreme_2006, coles_introduction_2001-1} for a detailed presentation), which models block maxima using the generalized extreme value (GEV) distribution. We propose an extreme quantile regression method where the GEV parameters $(\mu, \sigma, \xi)$ depend on covariates $x \in \mathcal{X}$. These parameters are estimated through a weighted version of the penalized maximum likelihood estimator, originally introduced by \cite{coles_likelihood-based_1999}, where the weights are derived from generalized random forests \cite{athey_generalized_2019}.  The generalized extreme value distribution is given by  
\begin{equation}\label{GEV}
G_{(\mu,\sigma,\xi)}(z)=
\begin{cases}
\exp\left(-(1+\xi \frac{z-\mu}{\sigma})_{+}^{-\frac{1}{\xi}}\right), &  \text{if } \xi \neq 0,\\[6pt]
\exp \left(-\exp \left(-\frac{z-\mu}{\sigma}\right)\right), & \text{if } \xi = 0,
\end{cases}
\end{equation}
defined on $\{z \in \mathbb{R}: 1+\xi \tfrac{z-\mu}{\sigma}>0\}$, where $\mu \in \mathbb{R}$ is the location parameter, $\sigma \in \mathbb{R}_+^*$ is the scale parameter, and $\xi \in \mathbb{R}$ is the extreme value index. Inspired by \cite{gnecco_extremal_2024}, our approach simultaneously addresses three challenges:  extrapolation in the tails through the GEV distribution,  complex quantile structures via adaptive similarity weights, and high-dimensional predictors space through the partitioning properties of grf. Finally, the penalty function employed in this work ensures the asymptotic optimality of the maximum likelihood estimator in large samples and improves estimation accuracy in small-sample settings.

This article is structured as follows. Section~\ref{ext_quant_backgr} introduces the fundamental concepts related to extreme quantile regression, the block maxima approach that we will use in our method, and the technique of generalized random forests. Section~\ref{setup_model} is dedicated to a detailed description of our methodology, with a particular focus on estimation and validation techniques suited to the tails of distributions. Finally, in Sections \ref{simul_study} and  \ref{real_data_study}, we evaluate the effectiveness of the proposed method through simulation studies and an empirical application using salary census data from the United States for the year 1980 \cite{angrist_quantile_2006}.
	
	\section{Extreme quantile regression}\label{ext_quant_backgr}	
The literature on extreme quantile estimation relies on the asymptotic results of extreme value theory, which allow for extrapolation beyond the range of observed data \cite{de_haan_extreme_2006, coles_introduction_2001-1}. In the case where no covariate is available, consider a sample \((Y_1, \dots, Y_n)\) consisting of independent and identically distributed realizations of a random variable \(Y\) with a cumulative distribution function \(F\). The goal is to estimate the quantile \(\mathcal{Q}(\tau) = F^{-1}(\tau)\) for an extreme probability level \(\tau \in (0,1)\).  When this level depends on the sample size \(n\), let \(\tau = \tau_n\). A probability level is considered extreme if \(\tau_n \to 1\) and \(n(1-\tau_n) \to d \geq 0\) as \(n \to +\infty\). In this case, the number of observations exceeding the quantile \(\mathcal{Q}(\tau_n)\) becomes limited as \(n\) increases, making empirical estimation particularly difficult. Estimating these extreme quantiles requires extrapolation beyond the range of available data \cite{de_haan_extreme_2006}, which necessitates robust asymptotic theory to make these extrapolations accurately and reliably. Let \(\tau_0\) denote the intermediate probability level from which the quantile order becomes extreme.  

One of the main results concerning extrapolation is the asymptotic theorem by \cite{fisher_limiting_1928} and \cite{gnedenko_sur_1943}, which establishes the limiting distribution of the maximum of a sequence of independent and identically distributed random variables. These authors demonstrated that if there exist sequences \(a_n > 0\) and \(b_n \in \mathbb{R}\) such that
\begin{equation}\label{dom_attr}
	F^n(b_n + a_n x)  \to G(x), \quad as\quad n  \to \infty,
	\end{equation}
 where $F^n$ denotes the $n$-th power of the cumulative distribution function of the random variable $Y,$  then 
  \(G\) is a non-degenerate distribution given by 
\[
G_{\xi}(x) = \exp\left(-(1 + \xi x)^{-\frac{1}{\xi}}\right),
\]
for all \(x\) such that \(1 + \xi x > 0\) and \(\xi \in \mathbb{R}\). The case \(\xi = 0\) corresponds to the limit of \(G_{\xi}(x)\) as \(\xi \to 0\), and is given by \(G_0(x) = \exp(-e^{-x})\).  A cumulative distribution function \(F\) is said to belong to the domain of attraction of an extreme value distribution, denoted \(F \in \mathcal{D}(G_{\xi})\), when there exist sequences \(a_n > 0\) and \(b_n \in \mathbb{R}\) such that the equality (\ref{dom_attr}) is satisfied. Depending on the sign of \(\xi\), three domains of attraction are distinguished: the Fréchet domain of attraction (when \(\xi > 0\)), the Gumbel domain of attraction (when \(\xi = 0\)), and the Weibull domain of attraction (when \(\xi < 0\)). The shape parameter \(\xi\) plays a key role in characterizing the behavior of the upper tail of the distribution. The Weibull class is characterized by a finite upper bound, while the Gumbel and Fréchet classes feature infinite tails with distinct rates of decay. More specifically, the Fréchet distribution, which has heavy tails, decays polynomially, while the Gumbel class decays exponentially. This difference in tail behavior is crucial, as it reflects structural differences in the modeling of extreme values. The unification of these three classes into a generalized extreme value distribution, whose explicit form is given by equation (\ref{GEV}), offers the major advantage of allowing the data to determine the most appropriate family, thanks to the estimation of the \(\xi\) parameter using inference methods. Thus, the choice of tail behavior is made without the need for prior decisions. Additionally, the uncertainty surrounding the estimation of \(\xi\) allows for an assessment of the relative relevance of the three types of distributions for their application to the available data \cite{de_paola_gev_2018}. In practice, special attention is given to the estimation of this parameter.

In this work, the GEV distribution is used to facilitate extrapolation in the tail of the distribution and thus the estimation of extreme quantiles. By inverting this distribution, the quantile corresponding to an extreme probability level \(\tau_n > \tau_0\) is obtained, expressed by
		\begin{equation}\label{GEV_noncondquant}
	Q(\tau) =
	\begin{cases}
	\mu + \dfrac{\sigma}{\xi} \left(\left(-\ln(\tau)\right)^{-\xi} - 1\right) & \text{if } \xi\neq 0, \\
	\mu - \sigma \ln\left(-\ln(\tau)\right) & \text{if } \xi = 0.
	\end{cases}
	\end{equation}
The estimation of the extreme quantile involves determining the parameters \(\mu\), \(\sigma\), and \(\xi\). Many works have focused on estimating these parameters while integrating statistical learning methods capable of addressing the issues raised in the introduction, namely the high dimensionality of the predictor space and the complexity of the relationships between the response variable and the explanatory variables. Within the framework of the Peak-over-Threshold (POT) approach in extreme value theory, \cite{velthoen_gradient_2023} proposed a method based on gradient boosting to estimate the parameters of the conditional generalized Pareto distribution (GPD), thereby providing an estimate of the extreme conditional quantile. Meanwhile, \cite{pasche_neural_2023} developed a method using neural networks to estimate these conditional quantiles. In a similar approach, \cite{gnecco_extremal_2024} introduced a method for estimating the parameters of the conditional GPD via a weighted version of the maximum likelihood estimator, with the weights determined by the generalized random forest method.
Although the maximum likelihood method is widely used and presents interesting asymptotic properties (\cite{bucher_maximum_2017, dombry_existence_2015, dombry_maximum_2019}), it is often criticized for its limited performance in the case of small samples. This weakness is a major obstacle in practice, where the analysis of extreme events often relies on a limited amount of data. Indeed, the rarity of extreme events means that, even over long periods of observation, the data available for fitting an extreme value model can remain scarce. To address this issue, \cite{coles_likelihood-based_1999} proposed a penalized maximum likelihood estimator for the parameters of the GEV distribution. This approach retains the asymptotic properties of the classical maximum likelihood estimator (MLE) while improving its robustness and performance in the presence of small samples. The likelihood function associated with the independent block maxima \( z_1, \dots, z_n \) is given by
\begin{equation}\label{MLE_equat}
L(\mu, \sigma, \xi) = \prod_{i=1}^n \frac{dG(z_i; \mu, \sigma, \xi)}{dz_i},
\end{equation}
and the maximum likelihood estimator  of the parameters \(\theta\) of the GEV distribution is given by  
$$\hat{\theta}=\arg\max\limits_{\theta \in \varTheta} L(\mu,\sigma,\xi),$$ with \(L(\mu, \sigma, \xi)\) given in equation (\ref{MLE_equat}) and \( \varTheta \subset \mathbb{R} \times \mathbb{R}^*_+ \times \mathbb{R}\). To address the limitations of the maximum likelihood method with small samples, our approach incorporates the penalized likelihood function proposed by \cite{coles_likelihood-based_1999}, as detailed in Section~\ref{section_PWL}. Specifically, we develop an estimation procedure for the parameters of the conditional GEV distribution based on a weighted extension of the estimator introduced by \cite{coles_likelihood-based_1999}, where the weights are derived from the generalized random forest method described below.

\subsection*{Generalized Random Forest}\label{sect_grf} 

The generalized random forest (grf), introduced by \cite{athey_generalized_2019}, extends the classical random forest method \cite{breiman_random_2001}. 
Like traditional random forests, grf aggregates predictions from $B$ decision trees, each grown on a bootstrap sample and built using random feature selection at each split, which promotes diversity and improves generalization.  The key novelty of grf is the flexibility to adapt the loss function guiding tree construction, making it suitable for a wide range of estimation tasks, such as conditional means, conditional quantiles \cite{athey_generalized_2019}, or extremal quantiles (see \cite{gnecco_extremal_2024} and \cite{Vidagbandji2025}).   Let $\{(X_i,Y_i)\}_{i=1,\cdots,n}$ be the training data. 
For a test point $x \in \mathcal{X}$, each tree provides an estimate of the conditional mean: $
\eta_b(x) = \sum_{i=1}^n \frac{\mathbb{1}_{\{X_i \in R_b(x)\}}}{|\{j : X_j \in R_b(x)\}|} Y_i,$  
 where $R_b(x)$ denotes the region (or leaf) of the $b$-th tree containing $x$. 
Equivalently, $R_b(x)$ corresponds to the set of training observations $X_i$ that fall in the same terminal node as $x$. 
The forest prediction $\eta(x)$ is obtained by taking the average of all the predictions from the $B$ trees  
\[
\eta(x) = \frac{1}{B} \sum_{b=1}^B \eta_b(x) 
= \sum_{i=1}^n w_n(x, X_i) Y_i,
\]  
with similarity weights defined as  
\begin{equation}\label{forest_weight}
w_n(x, X_i) = \frac{1}{B} \sum_{b=1}^B 
\frac{\mathbb{1}_{\{X_i \in R_b(x)\}}}{|\{j : X_j \in R_b(x)\}|},
\end{equation}
where $|E|$ denotes the cardinality of the set $E$.  The weights $w_n(x,X_i)$ quantify the contribution of observation $X_i$ to the prediction at $x$. 
They can be understood as the normalized frequency with which $X_i$ and $x$ fall in the same leaf across the ensemble of trees. 
Thus, $w_n(x,X_i)$ defines an adaptive similarity measure between $x$ and $X_i$, reflecting how strongly $X_i$ influences the prediction at $x$. 
This provides an intuitive interpretation: $R_b(x)$ identifies the local neighborhood of $x$ within a tree, and $w_n(x,X_i)$ aggregates these neighborhoods across the forest.   
 
In standard random forests, the similarity weights implicitly favor observations $X_i$ with $\mathbb{E}(Y|X=X_i) \approx \mathbb{E}(Y|X=x)$. 
However, this mechanism may fail to capture heterogeneity in quantile functions: an observation $X_i$ may have a low weight even if $\mathbb{Q}_{Y|X=X_i}(\tau) \approx \mathbb{Q}_{Y|X=x}(\tau)$ (see \cite{athey_generalized_2019}).  
By contrast, grf incorporates a loss function tailored to quantiles, ensuring that the similarity weights emphasize observations relevant for conditional quantile estimation.   Moreover, grf-based weights adaptively partition the covariate space in a way that highlights the most informative variables.  This is particularly advantageous in high-dimensional settings: unlike kernel or nearest-neighbor methods, which are severely affected by the curse of dimensionality, grf constructs data-driven neighborhoods around $x$ that capture complex and possibly nonlinear dependencies between the response and the covariates.  Therefore, weights based on the grf are a robust and flexible tool for approximating conditional quantiles, even in heterogeneous and high-dimensional environments.  They play a central role in the methodology we propose.  In what follows, these weights will be denoted by $w_i(x)$ for all $x \in \mathcal{X}$ and $i = 1, \dots, n.$

\section{Model and inference procedure}\label{setup_model}
In this section, we present in detail the procedure for estimating the parameters of the conditional GEV distribution, as well as our flexible method for estimating the conditional quantiles associated with extreme probability levels, in order to address the challenges discussed in section~\ref{ext_quant_backgr}.
\subsection{Setup for  extreme quantile regression}\label{setup_gev_Model}
Let \((X_1, Y_1), (X_2, Y_2), \ldots\) be a sequence of independent and identically distributed random vectors, coming from the random vector \((X, Y)\). Suppose that the conditional distribution function of \(Y\) given \(X=x\), denoted \(F_x(\cdot)\), belongs to the domain of attraction of the extreme value distribution, i.e., \(F_x \in \mathcal{D}\left(G_{\xi(x)}\right)\), with the corresponding normalizing sequences \(a_m(x)\) and \(b_m(x)\) defined in (\ref{dom_attr}). We divide the sequence of vectors into \(n\) blocks of size \(m\), such that the \(k\)-th block, for \(k \geq 1\), is given by
	\begin{equation*}
	B_{k,m} = \{(X_{(k-1)m+1}, Y_{(k-1)m+1}), \ldots, (X_{km}, Y_{km})\}.
	\end{equation*} 
Let \( Z_{k,m} = \max \{Y_i : (X_i, Y_i) \in B_{k,m}\} \) and \( X_{k,m} \) be the \( X_i \) associated with the \( Y_i \) maximizing \( \{Y_i : (X_i, Y_i) \in B_{k,m}\} \). Thus, for any \( m > 1 \), the sequence \( \{(X_{k,m}, Z_{k,m})\}_{k \geq 1} \) represents the block maxima, where the maximum is taken with respect to the \( Y \)-component. For any \( x \in \mathcal{X} \) and \( m \geq 1 \), the variables \( \{ Z_{k,m} | X_{k,m} = x \}_{k \geq 1} \) are independent and identically distributed with a distribution function \( F_x^m \). Thus, for any \( x \in \mathcal{X} \) and \( m > 1 \), we have
	\begin{equation*}
	\frac{Z_{k,m} - b_m(x)}{a_m(x)} \Big| X_{k,m} = x \xrightarrow{d} G_{\xi(x)}  \text{  as $m \to + \infty,$}
	\end{equation*}
with  
$
G_{\xi(x)}(z)= \exp\left(-(1+\xi(x)z)^{-\frac{1}{\xi(x)}}\right).
$ 
This allows us to approximate, for all \( x \in \mathcal{X} \), the distribution of \( Z_{k,m} | X_{k,m} = x \) by the conditional GEV distribution, where the parameters depend on the covariate \( x \). More specifically, the location, scale, and shape parameters are functions of \( x \), defined respectively by \( \mu(\cdot) : \mathcal{X} \to \mathbb{R} \), \( \sigma(\cdot) : \mathcal{X} \to \mathbb{R}_+^* \), and \( \xi(\cdot) : \mathcal{X} \to \mathbb{R} \). The conditional GEV distribution is given by
\begin{equation}\label{GEV_cond}
	G_{\theta(x)}(z)= \exp\left[-\left(1+\xi(x) \frac{z-\mu (x)}{\sigma (x)}\right)_+^{-\frac{1}{\xi(x)}}\right],
\end{equation}
with \( \theta(x) = (\mu(x), \sigma(x), \xi(x)) \) for all \( x \in \mathcal{X} \).  
The case \( \xi(x) = 0 \) corresponds to the limit of the expression given in (\ref{GEV_cond}) as \( \xi(x) \to 0 \), and is given by  
 \[
 G_{(\mu(x), \sigma(x), 0)}(z) = \exp\left[-\exp\left(- \frac{z - \mu(x)}{\sigma(x)}\right)\right].
 \] 
The quantile of order \( \tau \) (with \( \tau \) close to 1) of the conditional GEV distribution is obtained by inverting the distribution given in (\ref{GEV_cond}). It is given, for all \( x \in \mathcal{X} \), by  
	\begin{equation}\label{GEVquantile}
	Q_x(\tau) =
	\begin{cases}
	\mu(x) + \dfrac{\sigma(x)}{\xi(x)} \left(\left(-\ln(\tau^m)\right)^{-\xi(x)} - 1\right) & \text{if } \xi(x) \neq 0, \\
	\mu(x) - \sigma(x) \ln\left(-\ln(\tau^m)\right) & \text{if } \xi(x) = 0.
	\end{cases}
	\end{equation}
The conditional quantile depends on the three parameters of the conditional GEV distribution, namely \( \mu(x) \), \( \sigma(x) \), and \( \xi(x) \) for each \( x \in \mathcal{X} \). Therefore, its estimation requires the prior estimation of these parameters, which we propose to obtain initially, before substituting them into equation (\ref{GEVquantile}) to calculate an estimate of \( Q_x(\tau) \). There are various methods to estimate these parameters, among which the maximum likelihood estimator  is recognized for its flexibility and effectiveness in modeling extremes. However, this estimator performs worse than an alternative method based on probability-weighted moments (PWM) when applied to small samples. According to \cite{coles_likelihood-based_1999}, the superiority of the PWM method for small samples can be attributed to the assumption of a restricted parameter space, corresponding to finite population moments. In order to incorporate similar information into a likelihood-based approach, the authors propose a penalized maximum likelihood estimator, wich will be presented in section~\ref{section_PWL}. This estimator retains the flexibility of modeling and the asymptotic optimality of maximum likelihood, while improving its performance in the presence of small samples. The weighting procedure for this estimator, proposed in this work, which addresses the issues stated in the introduction, is detailed in the following section.	
 \subsection{ Penalized weighted likelihood estimator}\label{section_PWL}
As described in section~\ref{ext_quant_backgr}, the estimation of the conditional quantile presents two major challenges: the first is encountered when estimating quantiles for high probability levels (\(\tau\) close to 1), and the second arises when the quantile function is complex and the dimension of the predictor space is high. Our methodology addresses both of these challenges simultaneously. To facilitate extrapolation in the tail of the response variable \(Y\) conditional on \(X=x\), we rely on the asymptotic results of extreme value theory. More specifically, we model block maxima using the conditional GEV distribution, as defined in section~\ref{setup_gev_Model}. Furthermore, to capture the complex structure of the quantile function and ease the estimation in a high-dimensional covariate space, we use the weights \(w_i(x)\), which are estimated using the method of generalized random forests, detailed in section~\ref{sect_grf}. 
The limited performance of the maximum likelihood estimator  for estimating extreme quantiles in small samples can be attributed to the marked positive skewness in the distribution of the estimated parameter \( \xi \), which amplifies errors due to the nonlinear relationship between the quantile \( \mathcal{Q}_x(\tau) \) and \( \xi \), leading to significant biases. In contrast, the estimator based on probability-weighted moments (PWM) a priori assumes that \( \xi < 1 \), which limits the variation of \( \xi \) and generates a moderate negative bias, less severely penalizing errors (\cite{coles_likelihood-based_1999}). This bias-variance trade-off allows PWM to provide more reliable estimates of extreme quantiles.

	Since the maximum likelihood estimator may yield unstable results with small samples 
	and produce values of the shape parameter outside the admissible range, 
	\cite{coles_likelihood-based_1999} proposed the use of a penalty function. 
	This function incorporates into the likelihood the information that $\xi < 1$, 
	while making values close to 1 less probable than smaller values. 
	In this way, it imposes a coherent restriction on the parameter space, 
	analogous to that used for the PWM estimator, thereby improving the estimation 
	of GEV parameters in small-sample settings. The penalization acts as weak prior 
	information, ensuring consistent estimates without compromising the main advantages 
	of likelihood-based inference. The proposed penalty function is defined as follows
	\begin{equation}\label{funColes_pen}
	P_{\alpha, \lambda}(\xi) =
	\begin{cases}
	1 & \text{if } \xi \leq 0,\\
	\exp\left\{-\lambda\left(\tfrac{1}{1-\xi}-1\right)^{\alpha}\right\} & \text{if } 0<\xi<1,\\
	0 &  \text{if } \xi \geq 1,
	\end{cases}
	\end{equation}
	where $\alpha \geq 0$ and $\lambda \geq 0$ are tuning parameters. 
	Larger values of $\alpha$ impose a stronger relative penalty on values of $\xi$ close to 1, while $\lambda$ controls the overall weighting attached to the penalty.		
	To account for heterogeneous data structures, we propose an adaptive procedure to select the optimal values of the tuning parameters, based on the cross-validation method described in Appendix~\ref{VC_alphaLambda}. According to \cite{coles_likelihood-based_1999}, the combination $\lambda = \alpha = 1$ provides good performance for a wide range of $\xi$ values and sample sizes. Although our adaptive selection increases computational cost, it yields more stable and flexible performance in our applications, while also capturing heterogeneous structures in the data. 
	As illustrated in Figure \ref{fig_heatmap}, the optimal parameter values differ from the default values. However, it is important to note that, across the considered parameter grid, the variation in the cross-validation error is small, and therefore does not completely contradict \cite{coles_likelihood-based_1999} recommendation. These experiments thus confirm that the results remain robust to moderate changes in the tuning parameters. The resulting penalized maximum likelihood estimator (PMLE) provides more stable estimates of $\xi$ and, consequently, of extreme quantiles, enhancing the reliability of inference in extreme value analysis. For $\xi < 0$, the PMLE is practically indistinguishable from the classical MLE. In contrast, for $\xi \geq 0$, its behavior resembles that of the probability-weighted moments estimator, exhibiting lower variance at the cost of a slight negative bias (see \cite{coles_likelihood-based_1999}). Overall, in terms of the bias–variance trade-off, the PMLE performs at least as well as the probability-weighted moments estimator.

 The maximum penalized likelihood estimator retains all the asymptotic properties of the classical estimator, and that its performance in the presence of small samples is comparable to that of the PWM estimator. We incorporate this approach in the estimation of the parameters of the conditional GEV distribution in the context of this work. We define the estimator  \( \hat{\theta}(x) \) for the parameters of the conditional GEV, as the parameter \(\theta(x)\) that minimizes the weighted and penalized (negative) log-likelihood, given by
\begin{equation}\label{weigth_log_Coleslikelihood}
L^{pen}_n(\theta;x)=\sum_{i=1}^{n}w_i(x)\ell_{\theta}(z_i)+ \log(P_{\alpha,\lambda}(\xi)), \quad \forall x \in \mathcal{X},
\end{equation}
where \( \ell_{\theta}(z_i) = -\log\left( \frac{dG(z_i; \mu, \sigma, \xi)}{dz_i} \right) \) and \( P_{\alpha, \lambda}(\xi) \) is the penalty function  defined in (\ref{funColes_pen}).
 The weights \( w_i(x) \) for \( i = 1, \ldots, n \) are obtained through the grf method, which is described in section~\ref{sect_grf}. 
 The theoretical guarantees of the maximum likelihood estimator for the GEV distribution are well documented in the literature. The existence and convergence of the estimator are proved by \cite{dombry_existence_2015} and \cite{bucher_maximum_2017} for \( \xi > -1 \), and the asymptotic normality is demonstrated by \cite{dombry_maximum_2019} for \( \xi~>~ -\frac{1}{2} \).   Regarding the penalized estimator, \cite{coles_likelihood-based_1999} have shown that the penalized maximum likelihood retains the asymptotic properties of the classical MLE, since the penalty only acts as weak prior information when $\xi$ close to 1. In our framework, we extend this estimator to the weighted setting, where the weights are obtained from generalized random forests. Asymptotic guarantees for grf-based weights have been established by \cite{athey_generalized_2019}, who proved consistency and asymptotic normality for localized parameter estimates. Combining these results with the properties of the penalized MLE, we expect the weighted penalized estimator to inherit the same asymptotic guarantees under regularity conditions. A complete theoretical analysis is beyond the scope of this paper, but we acknowledge the importance of this point and plan to address it in future work.

 The likelihood function of the GEV distribution does not have a global maximum, but, by abuse of language, we will refer to \( (\hat{\mu}(x), \hat{\sigma}(x), \hat{\xi}(x)) \) as a weighted penalized maximum likelihood estimator if \( L^{\text{pen}}_n(\theta;x) \) has a local minimum at \( (\hat{\mu}(x), \hat{\sigma}(x), \hat{\xi}(x)) \). The parameter space of the GEV distribution, denoted \( \theta(\mathcal{X})=\{\zeta \in \mathbb{R} \times (0, +\infty) \times (-1, +\infty): \zeta =\theta(x)\text{ for some } x \in \mathcal{X} \} \), is unknown in practice. We define our estimator \( \hat{\theta}(x) \) as the value of the parameter \( \theta \) that optimizes (\ref{weigth_log_Coleslikelihood}) over a compact set \( \varTheta \subset \mathbb{R} \times (0, +\infty) \times (-1, +\infty) \) such that \( \theta(\mathcal{X}) \subset \varTheta \). It is thus defined by 
\begin{equation}\label{MLE_pen}
\hat{\theta}(x)\in \arg\min_{\theta \in \Theta} L^{pen}_n(\theta;x).
\end{equation}
These estimated parameters are then substituted into (\ref{GEVquantile}) and the final estimation of the conditional quantile of order \( \tau > \tau_0 \) is given by
 \begin{equation}\label{Q_esti_pen}
\hat{Q}_x(\tau) =
\begin{cases}
\hat{\mu}(x) + \dfrac{\hat{\sigma}(x)}{\hat{\xi}(x)} \left(\left(-\ln(\tau^m)\right)^{-\hat{\xi}(x)} - 1\right) & \text{if } \hat{\xi}(x) \neq 0, \\
\hat{\mu}(x) - \hat{\sigma}(x) \ln\left(-\ln(\tau^m)\right) & \text{if } \hat{\xi}(x) = 0.
\end{cases}
\end{equation}
Based on the work of \cite{gnecco_extremal_2024}, \cite{pasche_neural_2023}, and \cite{velthoen_gradient_2023}, we adopt in this study the value \( \tau_0 = 0.8 \), which corresponds to the probability level beyond which \( \mathcal{Q}_x(\tau) \) is considered an extreme quantile. Thus, our method applies to probability levels \( \tau \geq \tau_0 \). The algorithm outlining our approach is presented below in  algorithm~\ref{algo_GEV_erf}.

\subsection*{Algorithm}
Let \(\mathcal{D}_N = \{(X_i, Y_i)\}_{i=1 \cdots N}\) denote the initial sample, and \(\mathcal{D'}_n = \{(X_i, Z_i)\}_{i=1 \cdots n}\) the sample of block maxima defined in Section~\ref{ext_quant_backgr}, where \(m\) represents the block size, and \(\upsilon\) is the vector containing the hyperparameters associated with the generalized random forest.
\begin{algorithm}[h!]
	\caption{erf\_Pen}
	\label{algo_GEV_erf}
	\begin{algorithmic}[1]  \vspace{0.4cm}
		\BEGIN  $ $ \textbf{erf\_Pen-fit}
		\Input $ (\mathcal{D}_N$, $m$, $\upsilon$)  \vspace{0.1cm}
		\STATE   $\mathcal{D'}_n\leftarrow makeBloc(\mathcal{D}_N,m),$\vspace{0.1cm}
		\STATE $w_i(.)\leftarrow grf(\mathcal{D'}_n,\upsilon),$ \vspace{0.1cm}
		\Output  erf\_Pen $\leftarrow (\mathcal{D'}_n,w_i(.),m).$	
	\end{algorithmic}	\vspace{0.4cm}
	\begin{algorithmic}[1]
		\BEGIN $ $ \textbf{erf\_Pen-predict}\vspace{0.1cm}	
		\Input $ $ (erf\_Pen, x, $\tau$, $\alpha$, $\lambda$)\vspace{0.1cm}
		\STATE   $\hat{\theta}(x)\leftarrow \arg\min\limits_{\theta \in \Theta} L^{pen}_n{(\theta;x)}$ as defined in (\ref{MLE_pen}),
		\STATE 
		$\hat{Q}_{x}(\tau) \text{ is computed by substituting } \hat{\theta}(x) \text{ into equation } (\ref{Q_esti_pen}).$\vspace{0.1cm}
		\Output $\big(\hat{\theta}(x), \hat{Q}_{x}(\tau)\big).$  \vspace{0.4cm}			
	\end{algorithmic}			
\end{algorithm}
 The function \textbf{makeBloc} divides the initial sample into blocks of a given size \( m \) and returns the sample of block maxima.  The grf function is used to adjust the weights using the generalized random forests method proposed by \cite{athey_generalized_2019} and defined in section~\ref{ext_quant_backgr}.

 \section{ Simulation Study}\label{simul_study}
In this section, we conduct simulation studies to demonstrate the performance of the proposed method and its ability to address the issues outlined in the introduction. An independent sample of size \(N = 90,000\) is generated for the random vector \((X, Y)\), where the covariate \(X \in \mathbb{R}^p\) follows a uniform distribution over the hypercube \([-1, 1]^p\), and the conditional response variable \(Y | X = x\) follows a heavy-tailed distribution, according to the simulation study. The goal is to estimate the extreme conditional quantile \(\mathcal{Q}_x(\tau)\) for high probability levels, with \(\tau \in \{0.9, 0.99, 0.995, 0.9995,0.9999\}\). 
 The sensitivity analysis with respect to the block size is provided in Appendix~\ref{sect_sensitivity_analyseis}. This analysis indicates that the model is sensitive to the choice of block size $m$ and suggests an optimal range of $30-70$ for Scenario 1 and $35-110$ for Scenario 2, ensuring good agreement between the estimated GEV distribution and the block maxima. To ensure a sufficient number of observations for the learning process of our method, we set the block size to $m=40$, which produces $n=2,225$ observations from the $N=90,000$ generated data.

We compare our erf\_Pen method with two other widely used quantile regression methods in the literature, namely the quantile regression forest (qrf) method proposed by \cite{meinshausen_quantile_2006} and the generalized random forest (grf) method proposed by \cite{athey_generalized_2019}. We also consider a third method dedicated to extreme values, the Generalized Additive Extreme Value Model (evgam) introduced by \cite{Youngman2022}, in which the parameters of the GEV distribution are expressed as generalized additive functions. In our study, the location, scale, and shape parameters are modeled as smooth additive functions of the covariates, without including interaction effects. In order to demonstrate both the necessity and superiority of the penalized model, we included the unpenalized model (denoted by Unpen\_ERF) in our comparisons. This approach clearly illustrates the advantage of penalization. To evaluate the performance of the methods, we generate a test dataset \(\{x_i\}_{i=1, \dots, n'}\) with \(n' = 8,000\) (after forming the blocks with \(m = 40\), which gives 200 observations for the test) using the Halton sequence \cite{halton_algorithm_1964}. We use three metrics to compare the performance of the different methods. For each \(\tau \geq \tau_0\), we calculate the integrated squared error (ISE) for the estimated conditional quantiles \(\{\hat{ \mathcal{Q}}_{x_i}(\tau)\}_{i=1, \dots, n'}\) on the test data, as follows:
$
 ISE = \frac{1}{n'} \sum_{i=1}^{n'} \left(\hat{ \mathcal{Q}}_{x_i}(\tau) - \mathcal{Q}_{x_i}(\tau)\right)^2,
$
where \(x ~\longmapsto~\mathcal{Q}_{x}(\tau)\) represents the true quantile function for the probability level \(\tau\). By averaging over \(R = 100\) replicas of the fitting and estimation process for \(\mathcal{Q}_{x}(\tau)\), we obtain the mean integrated squared error (MISE), which is also used by \cite{gnecco_extremal_2024} and \cite{velthoen_gradient_2023} for evaluating the performance of their methods. The second metric we use is the integrated bias (IBias), introduced by \cite{wang_estimation_2013}, which is defined as the average: 
$
IBias = \frac{1}{n'} \sum_{i=1}^{n'} \left(\hat{ \mathcal{Q}}_{x_i}(\tau) - \mathcal{Q}_{x_i}(\tau)\right),
$
calculated on the test data. Finally, the third metric we use is the median absolute error (MedAE), which represents the median of the absolute differences \(\left|\hat{ \mathcal{Q}}_{x_i}(\tau) - \mathcal{Q}_{x_i}(\tau)\right|\) on the test data \(\{x_i\}_{i=1, \dots, n'}\).
 
The first scenario of this simulation study aims to demonstrate the ability of our method to provide accurate estimates when the predictor space is of high dimension and to test its robustness against noise. This scenario is similar to those presented by \cite{gnecco_extremal_2024} and \cite{velthoen_gradient_2023} in their simulation studies. The second scenario seeks to evaluate the performance of our method in complex situations, where the quantile function exhibits a highly nonlinear shape and where the covariate space is also large, thus addressing the issues raised in the introduction. In Appendix~\ref{Additional_Simul_study}, we include an additional simulation scenario to demonstrate the ability of erf\_Pen to provide accurate estimates when $\xi > 1$, which corresponds to a critical case in extreme value modeling. 
\begin{itemize}
	\item \textbf{Scenario 1:} We assume that \( Y|X=x \sim \gamma(x) \mathcal{T}_{\nu(x)} \), where \( \mathcal{T}_{\nu(x)} \) represents the Student's t-distribution with \( \nu(x) \) degrees of freedom. The dependence functions are defined by
	$$ \gamma(x) = 1 + \mathbb{1}_{x_1 > 0} \quad \text{and} \quad \nu(x) = 4 - (x_1^2 - x_2^2). $$ 
	In this scenario, only two variables, \( X_1 \) and \( X_2 \), are considered signal variables, while the remaining variables (\( p - 2 \)) are noise. This framework allows us to assess the performance of our method in a context where the predictor space is of high dimension, while also accounting for the presence of noise. The results obtained for this scenario are presented in section~\ref{sect_scenario1}.
	
	\item \textbf{Scenario 2:} We assume here that \( Y|X=x \sim \gamma(x) \mathcal{T}_{\nu(x)} \), where
	$$ \gamma(x) = 1 + 2\pi \varphi(2x_1, 2x_2) \quad \text{and} \quad \nu(x) = 3 + \frac{7}{1 + \exp(4x_1 + 1.2)}, $$  
	with \(\varphi\) represents the density of the centered bivariate normal distribution  
	with unit variance and a correlation coefficient of 0.75.
	In this scenario, we consider more complex forms for the functions \( \gamma(x) \) and \( \nu(x) \), which depend on the covariates \( x \). This scenario thus allows us to assess the robustness of our method against combined challenges, such as those mentioned in the introduction. The results of this analysis are presented in section~\ref{sect_scenario2}.
\end{itemize}
The performance of the generalized random forests method can be sensitive to the choice of hyperparameters, among which is the parameter \texttt{min.node.size}, which determines the minimum number of observations that a leaf can contain, and\texttt{ num.trees}. These hyperparameters, as well as the penalty parameters \( \lambda \) and \( \alpha \) in the penalized log-likelihood defined in (\ref{weigth_log_Coleslikelihood}), are optimized using the cross-validation method described in the appendix~\ref{VC_alphaLambda}.  In our simulation study,  we determined the optimal value of \texttt{min.node.size} from \( \{5, 10, 20, 50, 100\} \) using our cross-validation procedure. The other hyperparameters of the grf method were set to their default values in our simulation study. 
  
 \subsection{Scenario 1}\label{sect_scenario1}
 In this study, we  evaluate the performance of our method in a context where the predictor space is high-dimensional, and for different probability levels close to 1. To ensure a fair comparison between the methods, the learning and estimation process is carried out on the same data, consisting of the block maxima \(\{(X_i, Z_i)\}_{i=1,\dots,n}\).  Table ~\ref{table_MISE_DifMethod} illustrates the evolution of the logarithm of MISE as a function of the quantile level \( \tau \) for five conditional quantile estimation methods (erf\_Pen, Unpen\_ERF, evgam, grf, qrf) in Scenario $1.$ The results are presented for four configurations of covariate space dimensions: \( p \in \{10,20, 30, 40\}.\)

 \begin{table}[ht]
 	\centering
 	\begin{tabular}{|ccccc |c c|}
 		  \toprule
 		  \multicolumn{5}{|c|}{\textbf{Methods}} & &  \\
 		  \cmidrule(lr){1-5}
 		  erf\_Pen & Unpen\_ERF & evgam& qrf & grf & \textbf{$p$} & \textbf{$\tau$}  \\
 		  \midrule\hline
  		  
 		  3.544159 & 3.554158 & 3.559510 & 3.534113 & 3.528069&  10 & 0.8000 \\
  		  3.798552 & 3.827854 & 3.825582& 3.841730 & 3.831264 &  10 & 0.9000 \\
  		  4.726821 & 4.857544 &4.808876 & 5.049189 & 5.019799 &  10 & 0.9900 \\
  		  5.003168 & 5.175082 &5.104060 & 5.541602 & 5.513555 &  10 & 0.9950 \\
 		  5.020182 & 5.420251 &5.766500 & 6.166143 & 6.237731 &  10 & 0.9990 \\
 		  5.873965 & 6.223029 &6.040300 & 6.121009 & 6.219681 &  10 & 0.9995 \\\hline
 		  	   
 		3.476097 & 3.480686 & 3.455641 & 3.497357 & 3.506545 &  20 & 0.8000 \\
 		3.703584 & 3.717873 &3.689100 & 3.667943 & 3.669805 &  20 & 0.9000 \\
		4.522016 & 4.586443 &4.583590  & 4.651087 & 4.670408 &  20 & 0.9900 \\
 		4.760237 & 4.844819 &4.855621 & 4.874497 & 4.903185 &  20 & 0.9950 \\
		5.279417 & 5.421024 &5.466527 & 5.902284 & 5.895118 &  20 & 0.9990 \\		
 		5.487002 & 5.658167 &5.718320 & 6.355174 & 6.384959 &  20 & 0.9995 \\\hline

 		3.531561 & 3.541192 & 3.540931& 3.525474 & 3.529149 &  30 & 0.8000 \\
 		3.794515 & 3.821840 & 3.836681& 3.738508 & 3.741153 &  30 & 0.9000 \\
 		4.748053 & 4.860781 & 4.980860& 5.128974 & 5.143107 &  30 & 0.9900 \\
		5.031631 & 5.176482 &5.337251 & 5.571242 & 5.440019 &  30 & 0.9950 \\	
 		5.665040 & 5.895358 &6.159317 & 6.883641 & 6.927317 &  30 & 0.9990 \\	
 		 5.925686 & 6.197717 &6.508989 & 7.145184 & 7.259607 &  30 & 0.9995 \\\hline

 		3.469519 & 3.471978 &3.474885 & 3.443556 & 3.439969 &  40 & 0.8000 \\
		3.681014 & 3.688681 &3.698054 & 3.682495 & 3.673479 &  40 & 0.9000 \\
		4.432620 & 4.467006 &4.542034 & 4.472557 & 4.462648 &  40 & 0.9900 \\		
 		4.645826 & 4.690889 &4.793666 & 4.684846 & 4.685848 &  40 & 0.9950 \\		
 		5.095707 & 5.170992 &5.348638 & 5.598266 & 5.659914 &  40 & 0.9990 \\		
 		5.267733 & 5.358748 &5.572876 & 5.568182 & 5.640856 &  40 & 0.9995 \\\hline
 		 	\end{tabular}
 \caption{Log(MISE)  for different methods under varying dimensions $p $ and probability levels $\tau$.}
 \label{table_MISE_DifMethod}
 	 \end{table} 
The log(MISE) error increases gradually with \( \tau \), reflecting the growing difficulty in estimating extreme quantiles. This trend is particularly pronounced for \( \tau \geq 0.9 \), where the error rises almost exponentially, indicating increased variance and instability in predictions for these quantiles, as discussed in the introduction.
	For moderate probability levels (\( 0.8 \leq \tau \leq 0.9 \)), all methods exhibit similar performance, with nearly identical errors. In contrast, for higher quantiles (\( \tau > 0.9 \)), the erf\_Pen method stands out with lower error, especially as \( \tau \) approaches 1.
The evgam, grf, and qrf methods have slightly higher errors than the unpenalized Unpen\_ERF method, which  itself  is somewhat less accurate than
  erf\_Pen. For moderate dimensions $p$ of the covariate space, evgam and Unpen\_ERF have similar performance, but as $p$ increases, evgam becomes less accurate than Unpen\_ERF. Although evgam remains competitive compared to erf\_Pen across different probability levels, its accuracy decreases as $p$ increases. Furthermore, as pointed out by \cite{youngman_generalized_2019}, evgam becomes increasingly computationally demanding as 
	$p$ grows. These results indicate that our similarity-weighted estimation method provides better estimates than the other methods, and the integration of penalization further improves its performance.
	Increasing the dimension of covariates from \( p = 10 \) to \( p = 40 \) leads to a general increase in the log(MISE) error, confirming the challenge of estimating conditional quantiles in high dimensions. However, this increase remains relatively moderate for erf\_Pen, which retains an accuracy advantage, particularly for extreme values of \( \tau \).
	Thus, erf\_Pen appears to be the most robust method under these conditions, exhibiting a more controlled growth of error compared to Unpen\_ERF, evgam, grf, and qrf, making it a preferred choice for estimating conditional quantiles in high-dimensional settings and at extreme levels. This conclusion remains unchanged when performance is evaluated using other error metrics, as shown in the results presented in table~\ref{tabble_different_metrique}.

 \subsection{Scenario 2}\label{sect_scenario2}
\begin{figure}[h]
 	\centering
 	\includegraphics[width=0.9999\linewidth,height=0.6\textheight]{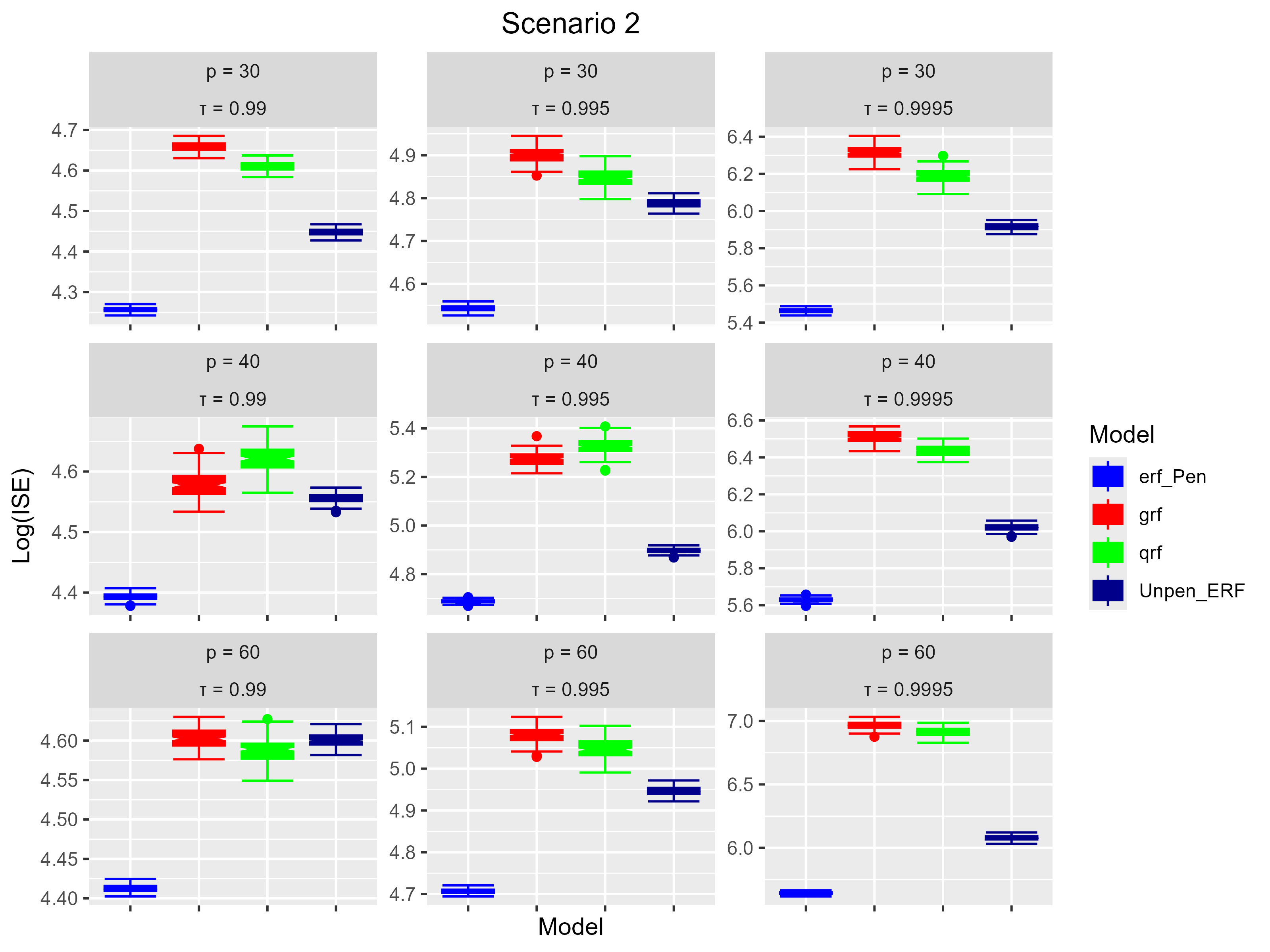}
 	\caption{Boxplots of Log(ISE) over 100 simulations for different values of \( p \) and extreme probability levels.}
 	\label{fig:plotlogisediffptau}
 \end{figure}

 In this scenario, we incorporate the unpenalize model (Unpen\_ERF) to compare its performance with proposed penalize method. Figure~\ref{fig:plotlogisediffptau} illustrates the comparison of logarithmic integrated squared errors (Log(ISE)) obtained by different conditional quantile estimation methods under Scenario 2 of our simulation study.
 It evaluates the robustness of the methods against extreme quantiles and the increase in the dimension \( p \) of the covariates. 
Due to the high computational cost of the evgam method when the covariate dimension $p$ is large, we exclude it from the comparison in this first analysis of this scenario. However, we include it in the second analysis summarized in Table~\ref{tabble_different_metrique}, for a fixed covariate dimension $p = 40$, to demonstrate that the results obtained in Scenario~1 remain valid for Scenario~2 across the different evaluation metrics considered. According to Figure~\ref{fig:plotlogisediffptau},
when \( \tau=0.99,\) the Unpen\_ERF, grf, and qrf methods yield almost similar performance, but the erf\_Pen method performs better, showing relatively moderate errors across the different covariate spaces considered, $p \in  \{30,40,60\}.$ 
 However, the Unpen\_ERF method, which does not incorporate penalization, although it shows better performance than the grf and qrf methods, exhibits larger errors than erf\_Pen, thereby confirming the usefulness of penalization for this type of estimation. 
  As \( \tau \) increases (\( \tau = 0.99, 0.995, 0.9995 \)), we observe an increase in error for all methods. However, the erf\_Pen method maintains more stable performance, with systematically lower error compared to grf and qrf, which are more affected by the extremity of the quantiles. This indicates that erf\_Pen is better suited for extreme quantile estimation. As the dimension of the covariate space increases (\( p =  30, 40,60 \)), the error tends to grow for all methods, reflecting the increased difficulty of estimation in high dimensions. However, erf\_Pen again stands out by showing a more controlled progression of error, while grf and qrf display more dispersed errors, indicating a loss of stability in their estimates.
 \begin{table}[h!]
 	\centering
 	\begin{tabular}{|c|c|c|c|c|c|}
 		\hline
 		\textbf{Method} & \textbf{MedAE} & \textbf{IBias} & \textbf{MISE} & \textbf{$\tau$} & \textbf{Scenario} \\ 
 		\hline
 		\textbf{erf\_Pen} & $6.20$             &  $6.31$             &  $3.68$             &             \multirow{7}{*}{\textbf{0.9}}  & \multirow{20}{*}{\textbf{ 1}} \\ \cmidrule{1-4}

 	\textbf{Unpen\_ERF}     & $ 6.22$              &  $ 6.33$      &  $ 3.69$             &              & \\ \cmidrule{1-4}
 	
 	\textbf{evgam}     & $6.20 $              &  $ 6.23$      &  $3.70 $             &              & \\ \cmidrule{1-4}
 	
 		\textbf{grf}     & $6.23$              &  $6.30$          &  $3.68$             &              & \\ \cmidrule{1-4}
 		\textbf{qrf}     &   $6.26$           &    $6.31$           &  $3.67$           &             & \\ \cmidrule{1-5}
 		\textbf{erf\_Pen} &     $9.87$         &     $10.1$          &    $4.64$           &            \multirow{7}{*}{\textbf{0.995}}  & \\ \cmidrule{1-4}
 		
 	\textbf{Unpen\_ERF}     & $10.1 $              &  $ 10.2 $      &  $ 4.69$             &              & \\ \cmidrule{1-4}
 	
 	\textbf{evagm}     & $10.1 $              &  $ 10.4$      &  $4.79 $             &              & \\ \cmidrule{1-4}

 		\textbf{grf}     &      $10.3$        &   $10.4$            &      $4.69$         &              & \\ \cmidrule{1-4}
 		\textbf{qrf}     &   $10.3$           &       $10.4$         &      $4.68$         &              & \\ 
 		\cmidrule{1-5}
 		
 		\textbf{erf\_Pen} &   $13.2$           &     $13.6$          &      $5.27$         &            \multirow{7}{*}{\textbf{0.9995}}  & \\ \cmidrule{1-4}

 	\textbf{Unpen\_ERF}     & $13.8 $              &  $14.2$             &  $5.36 $             &              & \\ \cmidrule{1-4}
 	
 	\textbf{evagm}     & $ 14.7$              &  $ 15.4$      &  $5.56 $             &              & \\ \cmidrule{1-4}

 		\textbf{grf}     &       $16.7$       &       $16.4$        &        $5.64$       &              & \\ \cmidrule{1-4}
 		\textbf{qrf}     &  $16.4$            &     $15.7$          &      $5.57$         &              & \\ 
 		\hline\hline

 		\textbf{erf\_Pen} &     $5.48$         &    $5.42$           &   $3.43$            &            \multirow{7}{*}{\textbf{0.9}}   & \multirow{20}{*}{\textbf{ 2}} \\ \cmidrule{1-4}

 		\textbf{Unpen\_ERF}     & $5.68 $              &  $ 5.60 $         &  $ 3.47$             &              & \\ \cmidrule{1-4}
 		
 		\textbf{evgam}     & $5.58 $              &  $ 5.57$      &  $ 3.46$             &              & \\ \cmidrule{1-4}

 		\textbf{grf}     &       $5.76$       &       $5.64$        &      $3.49$         &              & \\ \cmidrule{1-4}
 		\textbf{qrf}     &     $5.70$         &        $5.63$        &        $3.48$        &              & \\ \cmidrule{1-5}
 		\textbf{erf\_Pen} &      $9.06$         &    $9.13$            &     $4.63$           &            \multirow{7}{*}{\textbf{0.995}}  & \\ \cmidrule{1-4}

 		\textbf{Unpen\_ERF}     & $10.9 $              &  $ 10.9 $             &  $4.83 $             &              & \\ \cmidrule{1-4}

 		\textbf{evgam}     & $11.2 $              &  $11.5 $      &  $5.02 $             &              & \\ \cmidrule{1-4}

 		\textbf{grf}     &        $12.0$       &      $12.6$          &        $5.18$        &              & \\ \cmidrule{1-4}
 		\textbf{qrf}     &   $11.7$            &   $12.6$             &          $5.19$      &              & \\ 
 		\cmidrule{1-5}
 		\textbf{erf\_Pen} &       $11.8$        &        $11.7$        &    $5.49$            &            \multirow{7}{*}{\textbf{0.9995}}   & \\ \cmidrule{1-4}

 		\textbf{Unpen\_ERF}     & $15.6 $              &  $ 15.6 $             &  $5.80 $             &              & \\ \cmidrule{1-4}

 		\textbf{evgam}     & $16.7 $              &  $17.3 $      &  $5.82 $             &              & \\ \cmidrule{1-4}

 		\textbf{grf}     &    $17.2$           &      $17.8$          &      $ 5.91$          &              & \\ \cmidrule{1-4}
 		\textbf{qrf}     &         $17.0$      &     $17.4$           &      $5.84$          &              & \\ 
 		\hline
 		
 	\end{tabular}
 	\caption{Table of errors according to various metrics for each method and different scenarios.}
 	\label{tabble_different_metrique}
 \end{table}
In all configurations tested, erf\_Pen proves to be the most effective method. It exhibits significantly lower error than its competitors, particularly for extreme quantiles and in high dimensions. It is also more stable, with less variability in error compared to methods that show greater variability in their results. 
These results show that the erf\_Pen method is the most robust for extreme quantile estimation, even in the presence of numerous covariates. Its effectiveness is particularly notable for high values of \( \tau \), where it greatly outperforms  the grf method of \cite{athey_generalized_2019} and the qrf method of \cite{meinshausen_quantile_2006}. The increase in covariate dimension affects all methods, but erf\_Pen is better able to withstand this complexity. In summary, erf\_Pen stands out as the most suitable method for applications requiring reliable extreme quantile estimation, particularly in high-dimensional contexts. We also observe that our model remains effective in the presence of noise. 
The results presented in table~\ref{tabble_different_metrique} further confirm the superiority of the erf\_Pen method by using other error evaluation metrics. This table shows that in terms of MedAE and IBias, the erf\_Pen method outperforms the Unpen\_ERF, evgam, grf and qrf methods across different high probability levels, when the covariate dimension is set to \( p = 40 \). \textcolor{red}{}

 \section{Real dataset}\label{real_data_study}
We apply our methodology to predict the extreme quantiles of wage data from the 1980 U.S. Census. This dataset was used by \cite{angrist_quantile_2006} to illustrate a quantile regression method in non-extreme contexts, and more recently by \cite{gnecco_extremal_2024} in the analysis of extreme quantiles. The dataset consists of $65,023$ Black and White men, aged $40$ to $49$ years, with $5$ to $20$ years of education, and positive annual incomes and hours worked in the year prior to the census. In this study, we use weekly wages from $1979$ as the response variable \(Y\), expressed in U.S. dollars, calculated as the annual income divided by the number of weeks worked. The explanatory variables include a vector containing numerical variables representing age and years of education, as well as a categorical variable called Black, which is equal to $1$ if the respondent is Black and $0$ if the respondent is White.
 To increase the dimension of the predictor space and better assess the performance of our methodology in cases where the covariate space is large, we add $12$ random predictors, generated independently and uniformly within the interval \((-1, 1)\). This brings the total dimension of the predictor vector to \(p = 15\).
 \begin{figure}[h]
 	\centering
 	\includegraphics[width=0.99\linewidth, height=0.4\textheight]{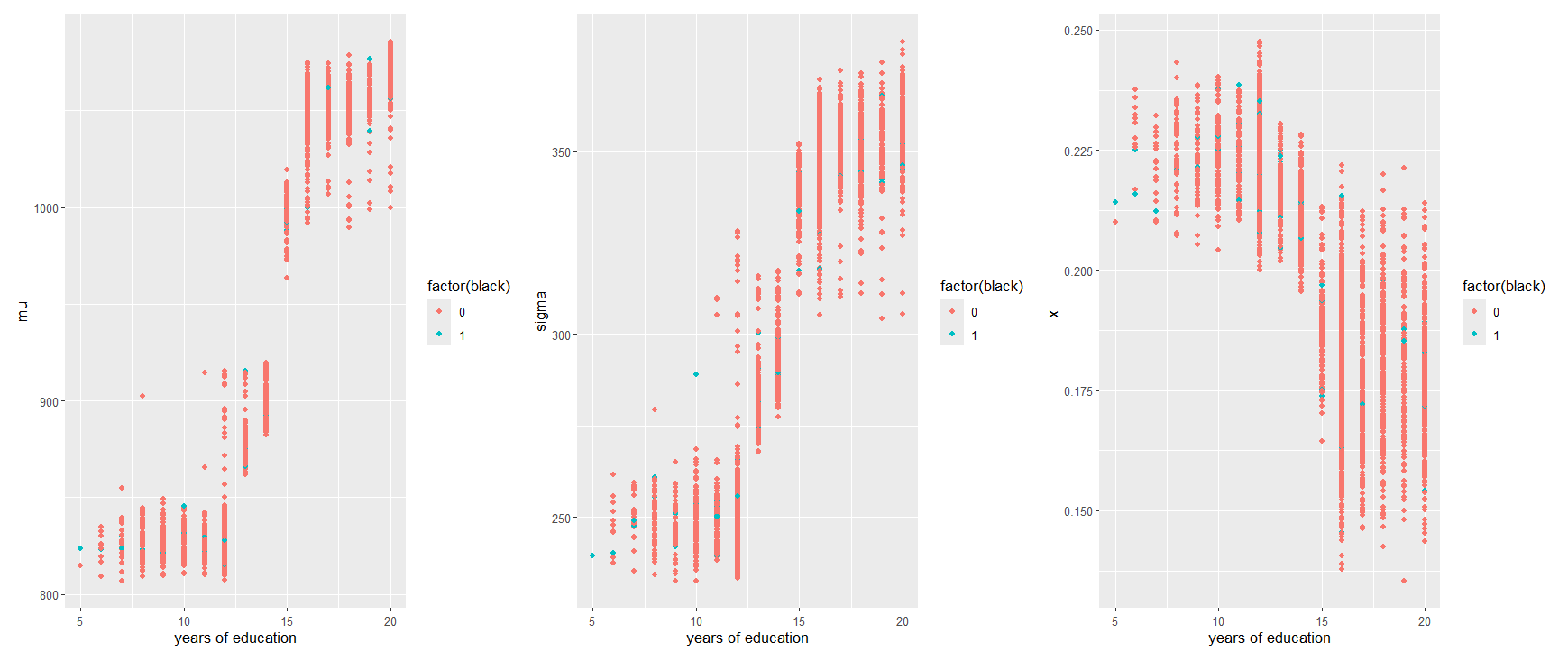}
 	\caption{Variations of the parameters \(\hat{\mu}(x)\), \(\hat{\sigma}(x)\), and \(\hat{\xi}(x)\) as a function of the number of years of education.}
 	\label{fig:plotxisigmgamyearedu}
 \end{figure}
  As done by \cite{gnecco_extremal_2024}, we split the dataset into two equal parts: the first part contains $32,511$ observations, used for exploratory analysis, and the other part contains $32,512$ observations, which is used for fitting and quantitatively evaluating the methods, including our proposed method, erf\_Pen, as well as the three other methods used in the simulation study (evgam from \cite{Youngman2022}, grf from \cite{athey_generalized_2019} and qrf from \cite{meinshausen_quantile_2006}).
 
For the exploratory analysis, we fit the erf\_Pen model to 40\% of the data, which corresponds to a total of $13,004$ observations, and estimate the parameters of the conditional generalized extreme value (GEV) distribution, \(\hat{\theta}(x) = (\hat{\mu}(x), \hat{\sigma}(x), \hat{\xi}(x))\), using the remaining 60\% of the data, applying the algorithm~\ref{algo_GEV_erf}. Since, as shown in Appendix \ref{fig_sensibility_analysis}, our method is sensitive to the block size, we used blocks of size $m=5$ observations in the exploratory analysis to improve estimation in the tail of the distribution. However, for the quantitative analysis, we applied the test statistic described in Appendix~\ref{sect_sensitivity_analyseis} and found that block sizes between 5 and 50 provide satisfactory goodness-of-fit. Therefore, we selected a block size of $m=20$ for the quantitative analysis. To determine the optimal values of the hyperparameters \(\alpha\) and \(\lambda\), which appear in the expression of the penalized likelihood given in (\ref{weigth_log_Coleslikelihood}), as well as \(min.node.size\) for the grf method, we apply cross-validation as described in appendix~\ref{VC_alphaLambda}. This procedure is performed for \((\alpha, \lambda) \in (0,5]\times(0,4]\) and \(min.node.size \in \{5,10,30,50\}\), keeping the default values for the other hyperparameters of the grf method. Finally, we proceed with the estimation of the parameters \(\theta\). Figure~\ref{fig:plotxisigmgamyearedu} shows the variation of the estimated parameters \(\hat{\mu}(x)\), \(\hat{\sigma}(x)\), and \(\hat{\xi}(x)\) as a function of the number of years of education. It is observed that \(\hat{\mu}(x)\) and \(\hat{\sigma}(x)\) increase with the level of education, while the shape parameter \(\hat{\xi}(x)\) follows an opposite trend, decreasing as the level of education increases. It is also noted that \(\hat{\xi}(x)\) takes positive values between 0.1 and 0.25, indicating heavy tails in the predictor space, thus confirming the analysis conducted by \cite{gnecco_extremal_2024}. Furthermore, the study highlights a balanced distribution of wages between Black and White Americans and shows that the parameters of the conditional GEV distribution are not influenced by age, as illustrated in figure~\ref{fig:plotxisigmgamage} in appendix~\ref{annex_dist_Para_age}.
 \begin{figure}[h!]
 	\centering
 	\includegraphics[width=0.99\linewidth, height=0.5\textheight]{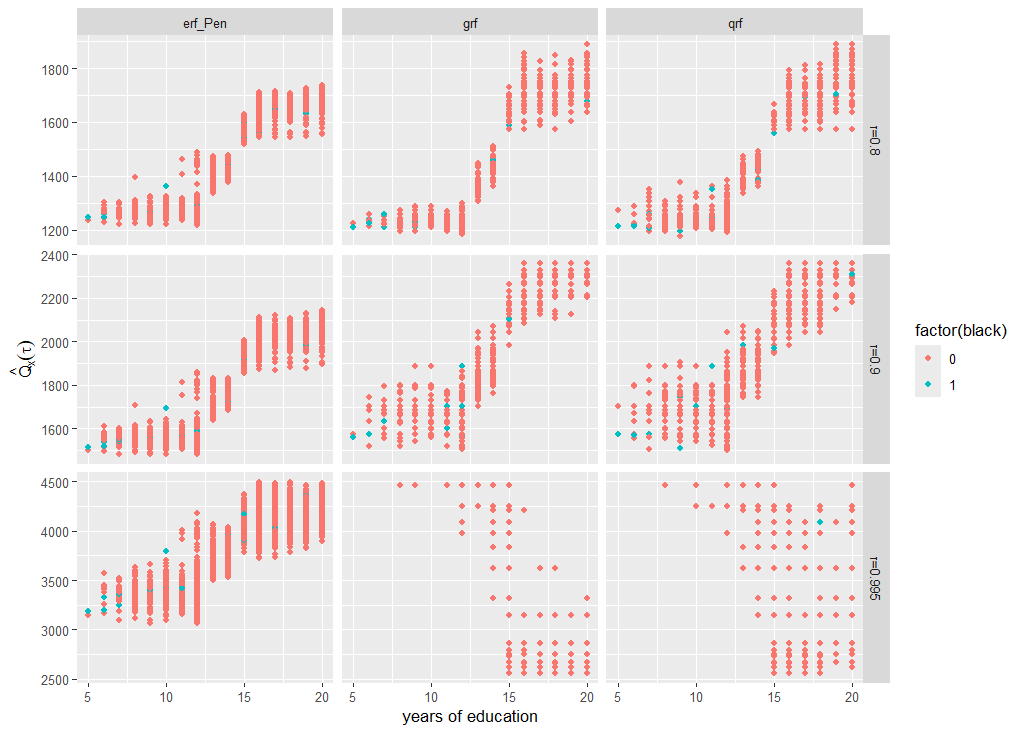}
 	\caption{Predicted conditionnal quantiles at level $\tau=0.8,0.9,0.995$ as fonction of years of eduction for erf\_Pen, grf and qrf method.}
 	\label{fig:plotestiquantileasdifftau}
 \end{figure} 
Figure~\ref{fig:plotestiquantileasdifftau} presents the estimates of the conditional quantiles of weekly income (\(\hat{Q}_{x}(\tau)\)) as a function of years of education, obtained by the three methods: erf\_Pen, grf, and qrf. The results are displayed for three quantile levels (\(\tau = 0.8\), \(\tau = 0.9\), and \(\tau = 0.995\)). The red and blue points correspond to individuals for whom the variable \texttt{black} is 0 or 1, respectively. An increase in income is observed with the number of years of education, a trend that becomes more pronounced at the highest quantiles, indicating that individuals with the highest income have spent significantly more years in education. When comparing the methods, our extrapolation method maintains a good shape of the quantile function even for high probability levels, which is not the case for the grf and qrf methods, where the shape deteriorates as \(\tau\) increases (particularly for \(\tau = 0.995\)). This highlights the ability of our method to capture the complex structure of the quantile function. This analysis underscores the growing impact of education on income.

After the exploratory analysis, we assess the quantitative performance of our method as well as the three other methods using the metric proposed by \cite{wang_estimation_2013}, which is also employed by \cite{gnecco_extremal_2024} to evaluate the performance of their own method. This metric is defined as follows:
\begin{equation}\label{colesMetrique}
	R_{n'}\left( \hat{Q}_x(\tau) \right) = \frac{\sum_{i=1}^{n'} \mathbb{1} \{ Y_i < \hat{Q}_{X_i}(\tau) \} - n'\tau}{\sqrt{n'\tau (1 - \tau)}}.
\end{equation}
The function $\hat{Q}_x(\tau)$ represents the estimated conditional quantile on the test sample $\{(x_i,y_i)\}_{i=1,\cdots,n'}$. 
	The metric $R_{n'}(\cdot)$ relies on the principle that, for a quantile at probability level $\tau,$ the proportion of observed values $y_i$  below the predicted quantile $\hat{Q}_{x_i}(\tau)$ should be close to $\tau.$		
	 Since $\mathbb{1}_{\{y_i < \hat{Q}_{x_i}(\tau)\}}$ has expectation $\tau$ and variance $\tau(1-\tau)$, this metric  is asymptotically normal according to the central limit theorem. A value of $R_{n'}(\cdot)$ close to zero indicates that the method provides predictions consistent with the theoretical definition of the quantile, while a large value reveals an underestimation or overestimation bias. This measure is particularly useful for the evaluation of extreme quantiles, as it captures the model’s ability to correctly reproduce the expected frequency of rare events, beyond standard criteria such as bias or mean squared error. For this purpose, we use the second part of the data, which was not used for the exploratory analysis, i.e., $32,512$ observations and consider the block size $m=20.$
 \begin{table}[h!]
 	\centering
 	\renewcommand{\arraystretch}{1.2}
 	\begin{tabular}{|c|c|c|c|c|}
 		\hline
 		\multirow{2}{*}{\centering \textbf{Model}} &\multicolumn{4}{c|}{\textbf{$\big|R_{n'}\left( \hat{Q}_x(\tau) \right)\big|$}} \\ 
 		\cmidrule{2-5}
 		
 					& $\tau = 0.8$  & $\tau = 0.9$ & $\tau = 0.995$  & $\tau=0.9999$ \\ \hline
 		\textbf{erf\_Pen}   &     3.951039  &  9.014595   &      59.755834 &    467.221579 \\ \hline
 		
 		\textbf{Unpen\_ERF}     &  4.205856     & 9.335311    &  60.925626   &   475.764389  \\ \hline

 		\textbf{evgam}     & 4.256064     &   9.200778   &  61.227295     & 483.950181    \\ \hline
 		
 		\textbf{grf}            & 3.210531       &  10.292044    &  68.164030     &   488.939327    \\ \hline
 		\textbf{qrf}            &  3.210713      & 10.302165    &   68.364223    &   488.282943  \\ \hline
 		
 	\end{tabular}
 	\caption{Performance of models based on the metrics defined in (\ref{colesMetrique}).}
 	\label{tab:performancequanti}
 \end{table}
We partition these data into 5 parts and, for each part $i \in \{1,2,3,4,5\}$, we train the model on $4/5$ of the observations, representing the set of observations excluding the $i$-th part. We then calculate the absolute value of $R_{n'}\left( \hat{Q}_x(\tau) \right)$ on the $1/5$ of the observations contained in the $i$-th part. After calculating this absolute error for each $i \in \{1,2,3,4,5\}$, we compare the average of these absolute errors for the different models and for different quantile orders.
 In this quantitative analysis, we included in our comparisons the unpenalized weighted likelihood estimator and evgam method. The table~\ref{tab:performancequanti} presents the mean values of these absolute errors  over the ten repetitions for different quantile levels. From this table, it is clear that our method provides good estimates compared to the methods of \cite{Youngman2022}, \cite{athey_generalized_2019} and \cite{meinshausen_quantile_2006}.    
The results show that the weighted likelihood estimator provides satisfactory estimates, whose accuracy improves when penalization is introduced. The evgam method remains competitive with the unpenalized version and, for high quantile levels ($\tau=0.995$ ), achieves performance comparable to that of erf\_Pen. Nevertheless, our approach remains stable at extreme quantile levels.

Although the real data application does not fall within the $\xi > 1$ regime, the additional simulation study presented in Appendix~\ref{Additional_Simul_study} highlights the limitations of the PMLE in the presence of very heavy tails. This is reflected in the larger values of $\log(\mathrm{MISE})$ compared to Scenarios~1 and~2. Figure~\ref{fig:additionnalsimulplot}, which reports the results for this critical setting, shows that our method provides the most accurate estimates across the different high-dimensional covariate configurations considered. Overall, these findings confirm the robustness and reliability of our approach across varying degrees of tail heaviness, including the challenging case $\xi>1$ examined in the simulation study.

\section{Conclusion}
In this paper, we propose a new flexible method, named erf\_Pen, for the reliable estimation of conditional quantiles at extreme probability levels. This method is also robust when the shape of the quantile function is complex and when predictors are potentially high-dimensional.  
Our methodology leverages the flexibility of generalized random forests to capture the complex structure of the quantile function, as well as the asymptotic results of extreme value theory \cite{fisher_limiting_1928, gnedenko_sur_1943} to facilitate extrapolation in the tail of the distribution via the conditional GEV distribution. The estimation of the latter’s parameters, using the penalty function of \cite{coles_likelihood-based_1999}, improves the estimation of the extreme value index, a key component in conditional quantile estimation.  Our approach stands out for its robustness and its ability to effectively capture the complex structure of the data while enabling reliable estimation in the tail of the distribution. Through a simulation study, we demonstrate that erf\_Pen outperforms other existing methods, notably the generalized random forests method of \cite{athey_generalized_2019}, the Generalized Additive Extreme Value Model of \cite{Youngman2022}  and the quantile regression forests method of \cite{meinshausen_quantile_2006}. It offers greater stability and improved accuracy for high quantile levels (\(\tau > 0.9\)), even in the presence of noise. Our results also show that erf\_Pen is more resistant to increasing covariate dimensionality, making it a particularly suitable alternative for high-dimensional applications.   An application to the 1980 U.S. Census wage data confirmed these empirical findings, demonstrating that our method preserves a good estimation of the quantile function, even for extreme quantiles (\(\tau = 0.995\)).  In conclusion, erf\_Pen represents a significant methodological advancement in the estimation of extreme conditional quantiles. Its performance and robustness make it a promising tool for applications in risk management, finance, economics, and other fields requiring reliable estimation of extreme conditional quantiles. Future work may include extending this approach to multivariate contexts as well as further theoretical analysis.
 
\backmatter
\section*{Acknowledgements}
We would like to thank the Editor, the Associate Editor, and the anonymous reviewers
for their helpful comments and suggestions, which improved this report of our work.

\section*{Declarations}

\textbf{Conflict of interest:} on behalf of all authors, the corresponding author states that there is no conflict of interest.

\begin{appendices}
	\section{Cross-validation method used to obtain \(\alpha\) and \(\lambda\) and the hyperparameters of grf.}\label{VC_alphaLambda}
	
	In this section, we present the method used to determine the tuning parameters of our model. Specifically, we address the selection of the penalization parameters \(\lambda\) and \(\alpha\), as well as the hyperparameters specific to the generalized random forest (such as \textit{min.node.size} and \textit{num.trees}). These parameters are chosen using a cross-validation procedure, adopting an approach different from the classical method described in section 7.10 of \cite{trevor_hastie_elements_2017}. We rely on the approach proposed by \cite{gnecco_extremal_2024}, which uses the likelihood function of the GPD as an error measure. However, in our case, we use the likelihood of the GEV distribution as the performance evaluation metric.
	
	To illustrate this method, let's consider a sample \(\mathcal{D}_n = \{(x_i, z_i)\}_{i=1,\cdots,n}\), available for our study. Suppose we wish to apply \(K\)-fold cross-validation. The first step is to partition the sample into \(K\) approximately equal-sized sub-samples, denoted \(\mathcal{D}^j\), with \(j = 1, \cdots, K\). A set of potential values is pre-defined for the penalization parameters \(\lambda\) and \(\alpha\), as well as for the hyperparameters specific to the generalized random forest (grf), such as \textit{min.node.size} and \textit{num.trees}. The cross-validation procedure then involves systematically evaluating all possible combinations of these parameters to identify those that optimize the model's performance. To do this, a search grid is constructed, denoted \(\{\varrho_1, \cdots, \varrho_S\}\), where each \(\varrho_s\) represents a vector containing a particular combination of the grf hyperparameters and the penalization parameters \(\lambda\) and \(\alpha\). For each parameter configuration \(\varrho_s \in \{\varrho_1, \cdots, \varrho_S\}\), the following steps are performed for each fold \(j \in \{1, \cdots, K\}\):
	\begin{enumerate}
		\item \textbf{Model training: } The model is trained on the \(K-1\) sub-samples, that is, on the entire dataset \(\mathcal{D}_n \setminus \mathcal{D}^j\), excluding the sub-sample \(\mathcal{D}^j\). This training is carried out using the \textit{Pen-erf-fit} algorithm, with \(\varrho_s\) as the configuration for the tuning parameters and the hyperparameters of the generalized random forest (grf).
		\item \textbf{Estimation of conditional parameters: }  the parameters \(\theta(x)\) associated with the conditional GEV distribution are estimated using the data from the sub-sample \(\mathcal{D}^j\).
		\item \textbf{Performance Evaluation:} A performance metric is calculated on the sub-sample \(\mathcal{D}^j\). In this context, the chosen metric is the negative log-likelihood of the GEV distribution.
	\end{enumerate}
	\begin{figure}[h!]
	\centering
	\includegraphics[width=0.99\linewidth, height=0.45\textheight]{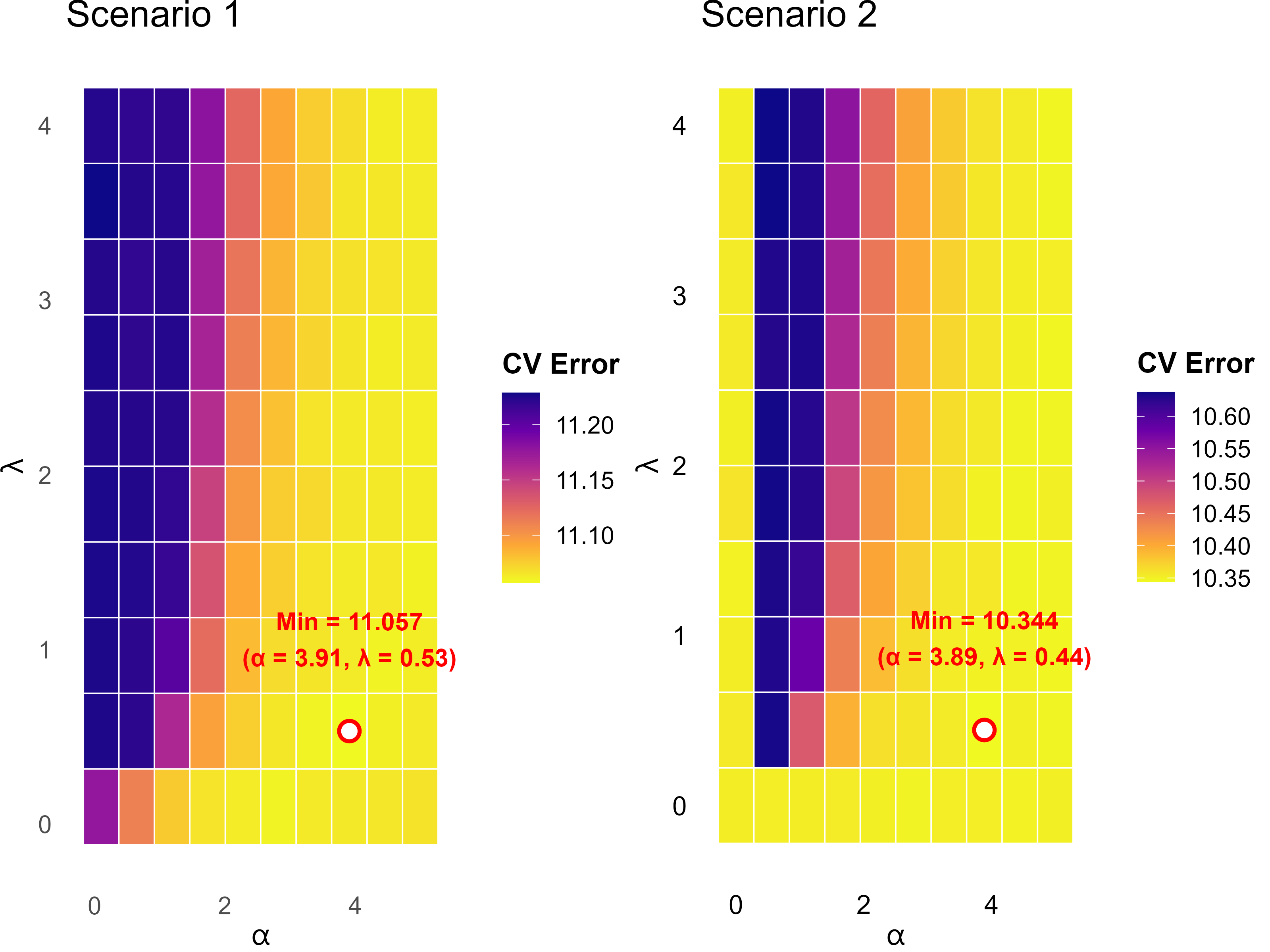}
	\caption{Cross-validation Error  across Scenarios as a function of $\alpha$ and $\lambda$ }
	\label{fig_heatmap}
\end{figure} 
	Once these steps are applied to the \(K\) sub-samples, the average error across the \(K\) partitions is calculated. This average error, called the cross-validation error (\textit{CV-error}), provides an overall estimate of the model's performance for a given configuration \(\varrho_s\). It is defined by:
	\[
	CV(\varrho_s) = \frac{1}{K} \sum_{j=1}^{K} \sum_{(x_i, z_i) \in \mathcal{D}^j} \ell_{(\hat{\theta}(x_i), \varrho_s)}(z_i),
	\]
	where \(\ell_{(\hat{\theta}(x), \varrho_s)}(z_i) = -\log\left(\frac{dG(z_i; \theta)}{dz_i}\right)\) represents the negative log-likelihood of the GEV distribution, with the parameters \(\hat{\theta}(x)\) estimated using \(\varrho_s\) as the parameter configuration. 
	This approach allows us to identify the optimal configurations of tuning parameters and hyperparameters by minimizing the cross-validation error. The optimal parameter \(\varrho_{optim}\) is the one that minimizes this error among all the configurations available in the set \(\{\varrho_1, \cdots, \varrho_S\}\). Formally, it is defined as: 
	\[
	\varrho_{optim} = \arg\min_{s \in \{1, \cdots, S\}} CV(\varrho_s).
	\]	
	This methodology ensures a robust and objective evaluation of the different parameter configurations, while identifying those that minimize the cross-validation error. This allows for reliable optimization of the model parameters, tailored to the specific characteristics of the data being studied.
		 The figure~\ref{fig_heatmap} illustrates how the cross-validation error varies with $\alpha \in (0, 5]$ and $\lambda \in (0, 4]$ across the two scenarios considered in this work. As shown in this figure, the optimal values of these parameters depend on the dataset used. However, the cross-validation error varies only moderately across the different parameter grids considered. This indicates that while satisfactory performance is obtained when $\alpha=\lambda=1,$ the proposed cross-validation procedure can still assist in selecting appropriate optimal values depending on the data.
		 		
	\section{Sensitivity analysis}\label{sect_sensitivity_analyseis}

In this section, we investigate the sensitivity of the erf\_Pen method with respect to the choice of block size. This parameter is crucial in the block maxima framework: a block size that is too small induces bias in the estimation, while a block size that is too large increases the variance of the estimator. The choice of $m$ must therefore achieve a bias–variance trade-off. To the best of our knowledge, there is no universally optimal procedure for selecting the block size, and practical applications often rely on natural choices such as monthly, quarterly, semi-annual, or annual blocks. 
\begin{figure}[h!]
	\centering
	\includegraphics[width=0.99\linewidth, height=0.5\textheight]{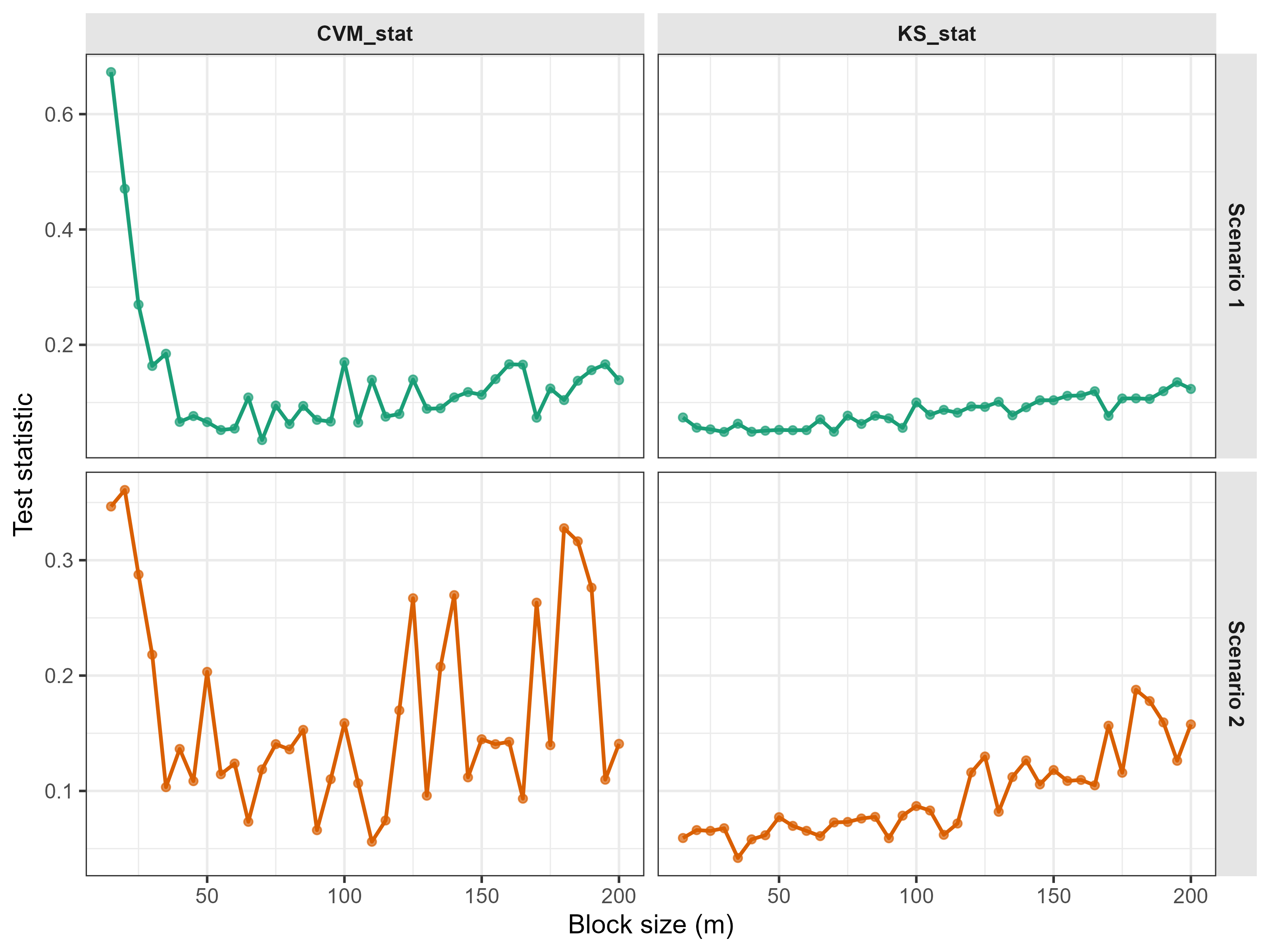}
	\caption{Evolution of test statistics as a function of block size under different scenarios  }
	\label{fig_sensibility_analysis}
\end{figure} 
We propose here an alternative approach for selecting the block size based on goodness-of-fit test statistics,  as done by \cite{wang_method_2016} in selecting block sizes for estimating extreme loads in engineering vehicles. Specifically, we fit the GEV model with different block sizes and evaluate the goodness-of-fit of the probability integral transform data against the uniform $(0,1)$ distribution using the Kolmogorov–Smirnov (KS) and Cramér–von Mises (CVM) tests.  The KS test measures the maximum deviation between the empirical and theoretical distributions, while the CVM test integrates the squared deviations over the entire support, thus providing a more global assessment. Their joint use yields complementary and robust evidence of goodness-of-fit, which is particularly useful for guiding the choice of block size.  

To this end, we consider a training sample $\mathcal{D}^{train}_N=\{(x_i,y_i)\}_{i=1}^N$  and an independent test sample $\mathcal{D}^{test}_l=\{(x_i,y_i)\}_{i=1}^l$.  We also consider a range of block sizes $m$ between 15 and 200, in increments of 5.  For each value of $m$, the model is fitted on the training sample  $\mathcal{D}^{train}_N$, and the parameters of the GEV distribution  $\left(\mu(x),\sigma(x),\xi(x)\right)$ are then estimated using the test sample $\mathcal{D}^{test}_l$.  Let $l^{\prime}$ denote the number of blocks of size $m$ formed from the test data. The probability integral transforms $\{G_{(\hat{\mu}(x_i),\hat{\sigma}(x_i),\hat{\xi}(x_i))}(z_i)\}_{i=1}^{l^{\prime}}$ 
should follow a uniform distribution if the estimation is adequate,  where $z_i$ denotes the maximum associated with the $i$-th block.   In this analysis, the penalization parameters are fixed at $\lambda=\alpha=1$, following the recommendations of \cite{coles_likelihood-based_1999},  and we adopt the default values of min.node.size and num.trees  for the generalized random forest method proposed by \cite{athey_generalized_2019}.
Figure~\ref{fig_sensibility_analysis} displays the evolution of the KS and CVM test statistics for the two scenarios considered. In Scenario~1 (top panels), the CVM statistic decreases sharply for block sizes smaller than $m=30$, then stabilizes at a relatively low and constant level for $m$ between 35 and 100. The KS statistic remains globally low but exhibits a slight upward trend as $m$ increases, suggesting a gradual loss of fit for very large block sizes.   In Scenario~2 (bottom panels), the results are more variable. The CVM statistic shows strong instability for small block sizes ($m < 35$), followed by marked fluctuations thereafter, indicating increased sensitivity of the estimation. The KS statistic exhibits a generally increasing trend with larger $m$, highlighting that excessively large block sizes deteriorate the goodness-of-fit.  

Overall, these  results confirm the method's sensitivity to block size and suggest an optimal range of $30-70$ for Scenario 1 and $35-110$ for Scenario 2, ensuring good agreement between the estimated GEV distribution and the block maxima. In practice, this type of sensitivity analysis can guide the empirical choice of $m$, although natural block sizes remain preferable. For our analyses, we selected $m=40$ in both scenarios, in order to balance model adequacy with a sufficient number of block maxima for the learning process.
\section{Additional Simulation Study}\label{Additional_Simul_study}

In this section, we conduct an additional simulation study to assess the ability of erf\_Pen to more accurately estimate conditional quantiles in a heavy-tailed regime where $\xi > 1$. This setting is particularly challenging, as some moments of the distribution may be infinite, making the estimation of extreme quantiles notoriously difficult. We consider a conditional distribution of the form
$
Y \mid X = x \sim \gamma(x)\,\mathcal{T}_{\nu(x)},
$
where $\mathcal{T}_{\nu(x)}$ denotes a Student's $t$ distribution with $\nu(x)$ degrees of freedom. As in Scenario~1, we set
$
\gamma(x) = 1 + \mathbb{1}_{\{x_1 > 0\}},
$
thereby introducing heterogeneity in the conditional scale. In this study, we assume a constant shape parameter given by
\[
\xi(x) = \frac{1}{\nu(x)} = 1.5, \quad \text{for all } x \in [-1,1]^p,
\] which corresponds to a heavy-tailed regime. This design allows us to evaluate the robustness of the proposed approach under a critical scenario in which extreme quantile estimation is particularly demanding. The covariates $X = (X_1, \dots, X_p)$ are generated independently from a uniform distribution on $[-1,1]^p$, with $p \in \{30,50\}$, consistent with the two scenarios considered in our simulation study. We evaluate the \texttt{log(MISE)} as a function of the probability level $\tau \in [0.8,1)$ for both covariate dimensions. We do not include the evgam model in this comparison because it becomes computationally expensive as $p$ increases. Moreover, although its performance is superior to that of the grf and qrf methods for moderate values of $p$, it becomes comparable to that of Unpen\_ERF, grf, and qrf for larger values of $p$, as evidenced by the results reported in Tables~\ref{table_MISE_DifMethod} and \ref{tabble_different_metrique}.

\begin{figure}[h!]
	\centering
	\includegraphics[width=0.99\linewidth, height=0.5\textheight]{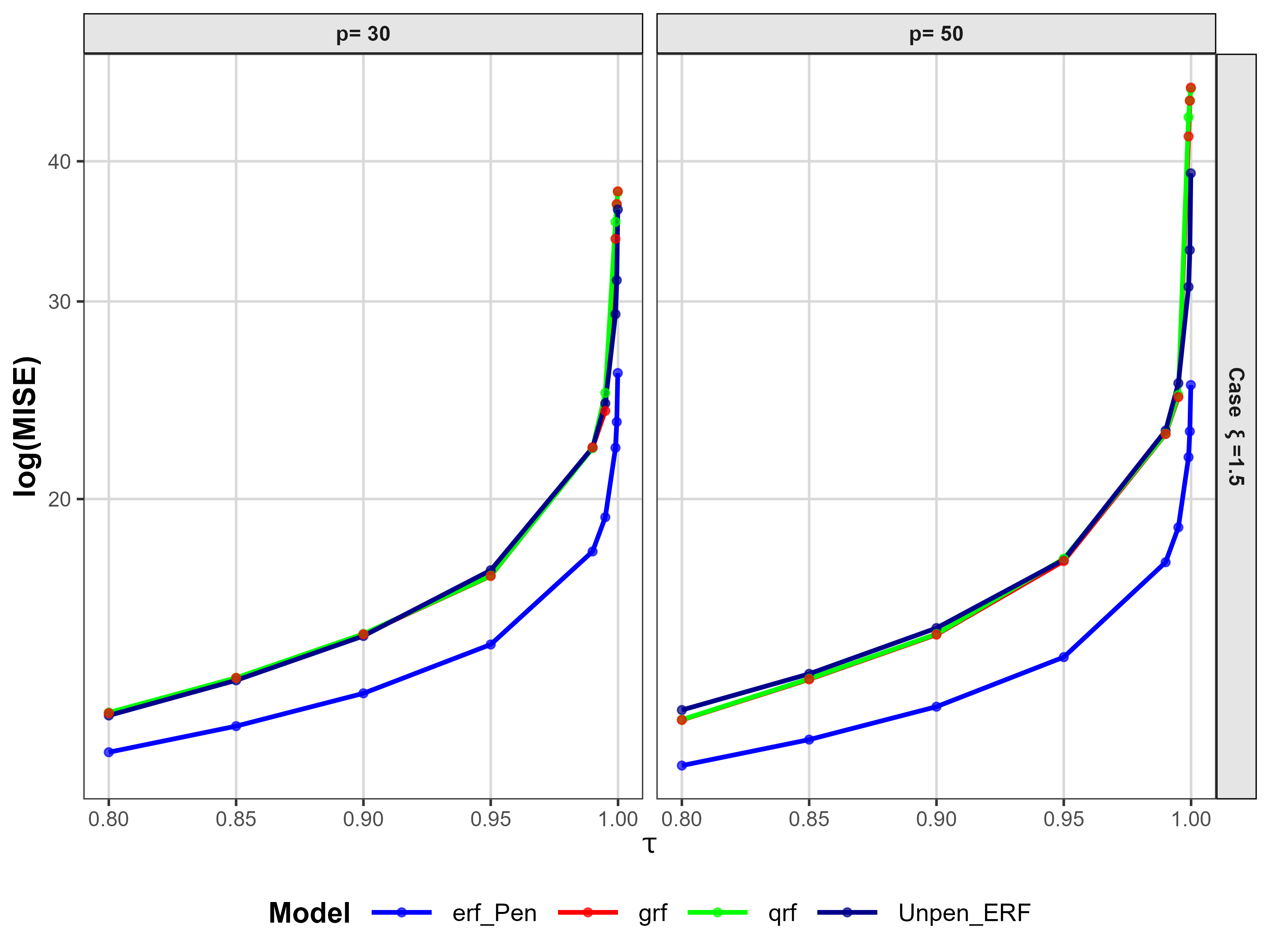}
\caption[Evolution of $\log(\text{MISE})$ as a function of $\tau$ for $p \in {20, 50}$]{Evolution of $\log(\text{MISE})$ as a function of $\tau$ for $p \in \{30, 50\}$}
\label{fig:additionnalsimulplot}
\end{figure}

The figure~\ref{fig:additionnalsimulplot} shows the evolution of log(MISE) as a function of probability levels $\tau$. The results shown in this figure indicate that erf\_Pen remains stable and delivers more accurate estimates than competing methods in this strongly heavy-tailed regime. In particular, the beneficial effect of penalization becomes more pronounced as both the probability level $\tau$ and the dimension $p$ increase, suggesting that the proposed approach is especially well suited for extreme quantile estimation when $\xi > 1$. We also observe that the non-penalized method performs similarly to qrf and grf for moderate dimensionality ($p = 30$), whereas for higher dimensionality ($p = 50$), it outperforms both competitors at high probability levels. These findings confirm that erf\_Pen remains applicable and effective beyond the standard setting $\xi < 1$, extending its relevance to heavy-tailed contexts with $\xi\geq1$.
	\section{Variation of the parameters \(\hat{\mu}(x)\), \(\hat{\sigma}(x)\), and \(\hat{\xi}(x)\) as a function of age. }\label{annex_dist_Para_age}	

	\begin{figure}[h!]
		\centering
		\includegraphics[width=0.99\linewidth, height=0.4\textheight]{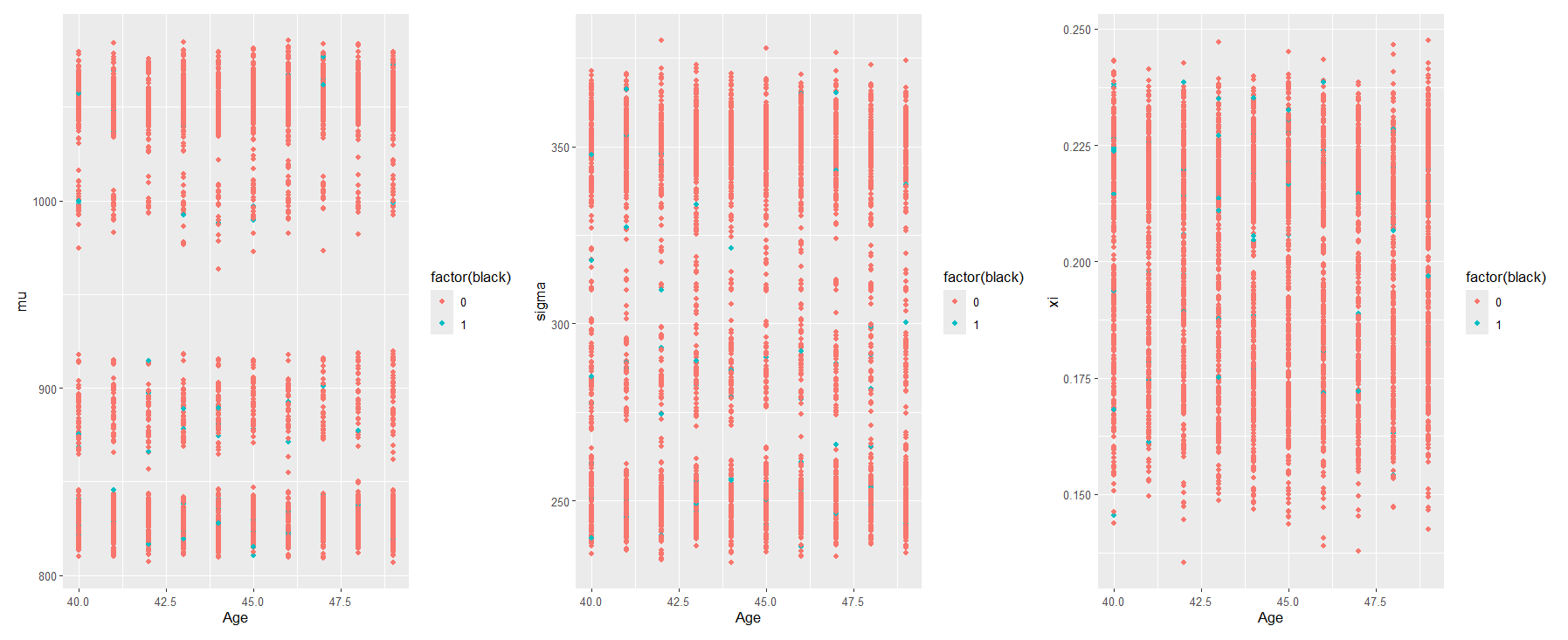}
		\caption{Variation of the parameters \(\hat{\mu}(x)\), \(\hat{\sigma}(x)\), and \(\hat{\xi}(x)\) as a function of age.}
		\label{fig:plotxisigmgamage}
	\end{figure}
	
\end{appendices}

 \bibliography{References}

 \end{document}